\documentclass[final]{siamart1116}
\usepackage{mathtools}
\usepackage{amsmath,amssymb,amsfonts,color}
\usepackage{algorithmic}
\usepackage{color}
\usepackage{booktabs}
\usepackage[mathscr]{eucal}
\usepackage{url}
\usepackage{dsfont}
\usepackage{tabularx}
\usepackage{bm}
\usepackage{xargs}

\newtheorem{remark}[theorem]{Remark}
\newtheorem{assumption}{Assumption}

\newcommand{\cA}{\mathcal{A}}
\newcommand{\cC}{\mathcal{C}}
\newcommand{\cF}{\mathcal{F}}
\newcommand{\R}{\mathbb{R}}

\newcommand{\cI}{\mathcal{I}}
\newcommand{\cG}{\mathcal{G}}
\newcommand{\cM}{\mathcal{M}}

\newcommand{\cU}{\mathcal{U}}

\newcommand{\cN}{\mathcal{N}}
\newcommand{\cJ}{\mathcal{J}}

\newcommand{\px}[1]{\partial_{{x} _{#1}}}
\newcommand{\pt}[1]{\partial_{t} #1}
\newcommand{\pxi}[1]{\partial_{\xi} #1}
\newcommand{\pdxi}[1]{\partial_{\xi\xi} #1}

\newcommand{\bU}{\mathbf{U}}
\newcommand{\bF}{\mathbf{F}}
\newcommand{\bG}{\mathbf{G}}
\newcommand{\bM}{\mathbf{M}}
\newcommand{\bb}{\mathbf{b}}
\newcommand{\bc}{\mathbf{c}}
\newcommand{\la}{\langle}
\newcommand{\ra}{\rangle}

\definecolor{philipp}{RGB}{205, 102, 000}

\title{Polynomial approximation of high-dimensional Hamilton-Jacobi-Bellman equations and applications to feedback control of semilinear parabolic PDEs}

\author{Dante Kalise\thanks{Imperial College London, Department of Mathematics, South Kensington Campus, London SW7 2AZ, United Kingdom ({\tt
dkaliseb@ic.ac.uk}).}
\and Karl Kunisch\thanks{University of Graz, Institute of Mathematics and Scientific
Computing, Heinrichstr. 36, A-8010 Graz, Austria and Johann Radon Institute for Computational and Applied Mathematics
(RICAM), Austrian Academy of Sciences, Altenberger Stra\ss{}e 69, 4040 Linz, Austria ({\tt
karl.kunisch@uni-graz.at}).
}}

\begin{document}
\maketitle

\begin{abstract}
	A procedure for the numerical approximation of high-dimensional Hamilton-Jacobi-Bellman (HJB) equations associated to optimal feedback control problems for semilinear parabolic equations is proposed. Its main ingredients are a pseudospectral collocation  approximation of the PDE dynamics, and an iterative method for the nonlinear HJB equation associated to the feedback synthesis. The latter is known as the Successive Galerkin Approximation. It can also be interpreted as  Newton iteration for the HJB equation. At every step, the associated linear Generalized HJB equation is approximated via a separable polynomial approximation ansatz. Stabilizing feedback controls are obtained from solutions to the HJB equations for systems of dimension up to fourteen.

\end{abstract}



\maketitle

%%%%%%%%%%%%%%%%%%%%%%%%%%%%%%%%%%%%%%%%%%%%%%%%%%%%%

\section{Introduction}\label{intro}

Optimal feedback controls for evolutionary control systems are of significant practical importance. Differently from open-loop optimal controls, they do not rely on knowledge of the initial condition and they can achieve design objectives, as for instance stabilisation, also in the presence of perturbations. Furthermore, the online synthesis of feedback control can be implemented in a real-time setting. It is well-known that their construction relies on special Hamilton-Jacobi-Bellman (HJB) equations, see for instance \cite{BCD97, FF14}. The solution of the HJB equation is the value function associated to the optimal control problem, and its gradient is used to construct  the optimal feedback control. In the very special, but important case of a linear control system with quadratic cost without constraints on the control or the state variables, the HJB equation reduces to a Riccati equation which has received a  tremendous amount of attention, both for the cases when the control system is related to ordinary or to partial differential equations. Otherwise one has to deal with the HJB equation which is a partial differential equation whose spatial dimension is that of the control system. Thus optimal feedback control for partial differential equations leads to HJB equations in infinite dimensions \cite{CLmult}.
After semi-discretization in space of the controlled partial differential equation (PDE), the HJB equation is posed in a space of dimension corresponding to the  spatial discretization  of the PDE \cite{F97}. For standard finite element or finite difference discretizations this leads to  high-dimensional HJB equations. This  is one of the instances which is referred to as {\em the curse of dimensionality} \cite{B61}.

Many attempts to tackle the difficulties posed for numerically solving the  HJB equations arising in optimal control have been made in the past or are currently being investigated. We refer, for instance, to \cite{FF14}, which mainly focuses on semi-Lagrangian schemes,  and further references given there.
A related approach to numerical optimal feedback control of PDEs is to semi-discretize the dynamics and to add a model order reduction step, either with Balanced Truncation or Proper Orthogonal Decomposition, in order to reduce the dimension of the dynamics to a number that is tractable for grid-based, semi-Lagrangian schemes. This approach has been successfully explored, for instance, in \cite{AF13,KK14,KVX04} and references therein. It strongly relies on a trustworthy representation of the dynamics via low-dimensional manifolds. Such a low-dimensional representation may deteriorate  when nonlinear and/or advection effects are relevant. Thus, it is important to strive for techniques, or combinations of techniques,  which allow to solve higher dimensional problems.

 Another direction of research  evolves around generalizing the Riccati-based approach to allow for nonlinearities in the state equation.  One such technique is termed state-dependent Riccati equation \cite{C97}. Here the coefficients in the 'ordinary' Riccati equation are functions of the state rather than constants as in the case of linear state equations. Another approach realizes the fact that the Riccati equation can be interpreted as the equation satisfied by the first term arising in the power series expansion of the value function, and attempts to improve by realizing also higher order terms in the expansion.  These methods are succinctly explained in \cite{BTB00}.

Yet another technique which has received a considerable amount of attention is termed  Successive Galerkin Approximation. Roughly speaking, the nonlinear HJB equation associated to the continuous-time optimal control problem is solved by means of a Newton method. At each iteration, the control law is fixed. This leads to a Generalized Hamilton-Jacobi equation (GHJB) which is linear. The iteration is closed by an update of the control law based on  the gradient of the value function. This method  was intensively investigated in \cite{BST7, BST98}, see also \cite{BTB00}, and the references given in these citations. It is worth mentioning that the discrete-time counterpart of this method corresponds to the well-known policy iteration or Howards' algorithm \cite{H60,BMZ09,AFK15}.

The numerical examples in \cite{BST7, BST98,BTB00} do not go beyond dimension five, and most, if not all, of the published numerical results for nonlinear HJB equations do not exceed dimension  eight \cite{BGGK13,GK16,G16}. An alternative sparse grid approach for high-dimensional approximation of HJB equations based on open-loop optimal control has been presented in \cite{KW16}, with tests up to dimension six. Numerical methods relying on tensor calculus have been shown to perform well in high-dimensional settings where the associated HJB equation is a linear PDE \cite{YPL}. { A key feature of these works is the use of sparse tensor products either for the construction of the basis or for the representation of the solution. This idea constitutes a cornerstone of high-dimensional approximation \cite{BM05,CCS14}, and their applicability ranges from sparse grid approximations \cite{BG04},  to polynomial chaos expansion \cite{SG11,HS13}  and uncertainty quantification \cite{GWZ17}.}

In the present paper, to solve optimal control problems for certain classes of semilinear parabolic equations we shall proceed as follows. To accommodate the curse of dimensionality, the discretization of the PDE is based on a pseudospectral collocation method, allowing a higher degree of accuracy with relatively few collocation points. To solve the resulting  HJB we utilize a Newton method based  on the GHJB equation as described above. Next, the discretization of the GHJB equation is addressed through a Galerkin approximation with polynomial, globally supported, ansatz functions. While this mitigates the curse of dimensionality in terms of removing the mesh structure, it leads to high-dimensional integrals. We therefore resort to separable representations for the system dynamics and for the basis set of the polynomial approximation. The separability assumption reduces the computation  of the Galerkin residual equation to products of one-dimensional integrals. The combination of these procedures allowed us to solve HJB equations related to nonlinear control systems up to dimension fourteen by means of basic parallelization tools. The successful use of the Newton procedure requires to provide a feasibly initialization, i.e. a sub-optimal, stabilizing control. Since we do not consider constraints, this is not restrictive for finite horizon problem, but can be challenging for infinite horizon problems, and specifically for the stabilization problems which are considered in the present paper. In this respect we developed a continuation procedure based on the use of  a discount factor. Specifically, we consider a nested iterative procedure: within the outer loop the value of a positive discount factor is driven to zero, within the inner loop the HJB equation is solved approximately  for a fixed discount factor.  With this approach, which, is summarized in Algorithms \ref{alg:sga1} and  \ref{alg:sga2} below, we managed to solve optimal feedback stabilization problems for semilinear parabolic equations with different stability  behavior of the desired steady state.

Let us give a brief outline of the paper. Section 2 sets the stage and provides the discussion of a special case to facilitate the understanding of the following material. In Section 3 the solution process of the HJB equation is detailed. In Section 4 we provide the formulas which are needed to numerically realize the discretized HJB equation after a separable basis has been chosen. Numerical experiments are documented in Section 5. There we can also find comparisons to suboptimal feedback strategies based on Riccati and asymptotic expansion techniques.

\section{Infinite horizon optimal feedback control}
We consider the following undiscounted infinite horizon optimal control problem:

\[\underset{u(\cdot)\in\cU}{\min}\;\cJ(u(\cdot),x_0):=\int\limits_0^\infty \ell(x(t))+\gamma|u(t)|^2\, dt
\]
subject to the nonlinear dynamical constraint
\[\dot x(t)= f(x(t))+g(x)u(t)\,,\quad x(0)=x_0,\] where we denote the state $x(t)=(x_1(t),\ldots,x_d(t))^t\in\R^d$, the control $u(\cdot)\in\cU$, with $\cU=\{u(t):\, \R_+\rightarrow U\subset\R^m\}$, the state running cost $\ell(x)>0$, and  the control penalization $\gamma>0$. Furthermore, we assume the running cost and the system dynamics $f(x):\R^d\rightarrow\R^d$ and {$g(x):\R^d\rightarrow\R^{d\times m}$} to be $\cC^1(\R^d)$. Throughout it is assumed that $f(0)=0$ and $\ell(0)=0$. Our focus is therefore asymptotic stabilization to the origin.

It is well-known that the optimal value function
\[V(x_0)=\underset{u(\cdot)\in\cU}{\inf} J(u(\cdot),x_0)\]
characterizing the solution of this infinite horizon control problem is the unique viscosity solution of the Hamilton-Jacobi-Bellman equation
\begin{equation}\label{hjb}
\underset{u\in U}{min}\{ DV(x)^t(f(x)+g(x) u)+ \ell(x)+\gamma|u|^2\}=0\,,\quad V(0)=0\,,
\end{equation}
with $D V(x)=(\px{1}V,\ldots,\px{d}V)^t$. Here we follow the convention of dropping the subscript of $x_0$. We study this equation in the unconstrained case, i.e., $U\equiv \R^m$, where  the explicit minimizer $u^*$ of \eqref{hjb} is given by
\begin{equation}\label{optc}
u^*(x)=\underset{u\in U}{argmin}\{ DV(x)^t(f(x)+g u)+ \ell(x)+\gamma|u|^2\}=-\frac{1}{2\gamma} g(x)^t DV(x)\,.
\end{equation}
note that by inserting this expression for the optimal control in \eqref{hjb}, we obtain the equivalent HJB equation
\begin{equation}\label{hjb2}
DV(x)^t f(x)-\frac{1}{4\gamma}DV(x)^tg(x)g(x)^tDV(x)+\ell(x)=0\,,
\end{equation}
which under further assumptions can be simplified to the Riccati equation associated to linear-quadratic infinite horizon optimal feedback control.

The methodology we present in this work is applicable to systems fitting the aforedescribed setting, although for the sake of simplicity we restrict the presentation by the following choices:
\begin{itemize}
\item[(i)]  the control $u(t)$ is a scalar variable, i.e. $m=1$.
\item[(ii)] the running cost $\ell(x)$ is quadratic, i.e. $ x^T Q x$, with $Q$ positive-definite,
\item[(iii)] the control term $g(x)\equiv g$ is a constant vector in $\R^d$.  %({\color{blue} can be extended to a separable structure in $x$ as for $f$}).
\end{itemize}
At this point, our setting differs from the linear-quadratic case as it allows nonlinear dynamics.  For the numerical scheme that we develop the following assumption is crucial:

\begin{assumption}\label{as:1}
The  free dynamics $f(x):\R^d\to\R^d, f(x):=(f_1(x),\ldots,f_d(x))^t$ { is a sum of separable functions} in every coordinate $f_i(x)$
\[f_i(x)=\sum_{j=1}^{n_{f}}\prod_{k=1}^{d}\cF_{(i,j,k)}(x_k)\,,\]
where $\cF(x):\R^d\rightarrow\R^{d\times n_f\times d}$ is a tensor-valued  function. In the case $g=g(x)$, then we shall also assume a similar separable structure for $g(x)$.
\end{assumption}

{ Separated representations are a fundamental tool for mitigating the curse of dimensionality, often leading to algorithms that scale linearly in $d$. Its computational efficiency depends on the number of summands or \textsl{separation rank} ($n_f$ above). In this work the we shall assume the separated representation is exact. This can be readily checked as we consider dynamical systems arising from the application of the method of lines to semilinear parabolic PDEs. This translates into working over an ODE control system of the form
$\dot x(t)=Ax(t)+N(x(t))+u(t)\,,$ where $A$ is a linear operator, and $N(x(t))$ is a polynomial source term. Even though such a setting is quite general and covers a wide class of meaningful problems, it excludes an important set of agent-based control models where the governing dynamics depend on a metric interaction between states, i.e. $f(x)=f(\|x_i-x_j\|)$ (see \cite{BFK15} and references therein for control-related examples). Nonetheless, it is possible to address the problem of finding an approximate best separated representation of fixed rank $n_f$, but this procedure and its error analysis is beyond the scope of the present work. The interested reader can find in \cite{BM05} a thorough presentation of this topic, with a concrete application to linear HJB equations in \cite{YPL}.}

{ Under the framework provided by Assumption \ref{as:1}, the methodology can be directly applied to multidimensional control signals, non-quadratic state costs, and state dependent $g(x)$ which corresponds to bilinear control systems. The additional computational cost of addressing multidimensional control signals scales linearly with the dimension of $u$, whereas the computational burden associated to a bilinear control system will depend on the separability degree of $g(x)$. Non-quadratic, separable state costs can incorporated at a negligible computational cost with straightforward modifications of our setting.}

\subsection{Towards optimal feedback control of semilinear parabolic equations}\label{pdesub} In the following, we illustrate how the presented framework sets  the grounds for a computational approach for approximate optimal feedback controllers for nonlinear PDEs. We consider the following optimal stabilization problem:

\begin{equation}\label{costex1}
\underset{u(\cdot)\in L^2([0;+\infty))}{\min}\;\cJ(u(\cdot,X_0):=\int\limits_0^\infty \|X(\cdot,t)\|^2_{L^2(\cI)}+|u(t)|^2\, dt
\end{equation}
subject to the semilinear parabolic equation
\begin{align}\label{ex1}
\pt X(\xi,t)&=\pdxi X(\xi,t)-X(\xi,t)^3+\chi_{\omega}(\xi)u(t)\,,\quad\xi\in\cI=[-1,1]\,, t\in\R^+,\\
\pxi X(-1,t)&=\pxi X(1,t)=0\,,\quad X(\xi,0)=X_0\,.\notag
\end{align}
In this case, the scalar control acts through the indicator function $\chi_{\omega}(\xi)$, with $\omega\subset\cI$. At the abstract level, this corresponds to an infinite-dimensional optimal control problem. A first step towards the application of the proposed framework is the space discretization of the system dynamics, leading to finite-dimensional state space representation. The use of the pseudospectral collocation methods for parabolic equations has been studied in \cite{os06,qvbook}, and leads to a state space representation of the form
\begin{equation*}
\dot X(t) = A X(t)-X(t)^3+B u(t)\,,
\end{equation*}
where the discrete state $X(t)=(X_1(t),\ldots,X_d(t))^t\in\R^d$ corresponds to the approximation of  $X(\xi,t)$ at $d$ collocation points $\xi_i= -cos(\pi i/d)$, $i=1,\ldots,d,$ and $X^3$ is the coordinatewise power. The matrices $A\in\cM^{d\times d}$ and $B\in \R^d$ are finite-dimensional approximations of the Laplacian and control operators, respectively. { Such a discretization of the dynamics directly fulfills the separability required in Assumption \ref{as:1}}, as the i-th equation of the dynamics reads
\[\dot X_i(t)=A_{i,1}X_1(t)+\ldots+A_{i,d}X_d(t)-X_i(t)^3+B_iu(t)\,,\]
with a separability degree $n_f=d+1$.
It is very important to note that semidiscretization in space of a wide class of time-dependent PDEs will lead to finite-dimensional state space representations of this type, thus the applicability of the presented framework is only limited by the dimensionality of the associated HJB equation. This motivates the choice of a pseudospectral collocation method for the discretization, as it is possible to obtain a meaningful representation of the dynamics with considerably fewer degrees of freedom than classical low-order schemes. However, if pseudospectral collocation is not a suitable discretization method for the dynamics, model reduction procedures such as balanced truncation, proper orthogonal decomposition, or reduced basis techniques shall also lead to separable state-space representations. Once the finite-dimensional state state space representation is obtained, we proceed to approximate the solution of the associated HJB equation \eqref{hjb}, leading to the optimal feedback controller \eqref{optc}.

\begin{figure}[!ht]
\centering
\includegraphics[width=0.45\textwidth]{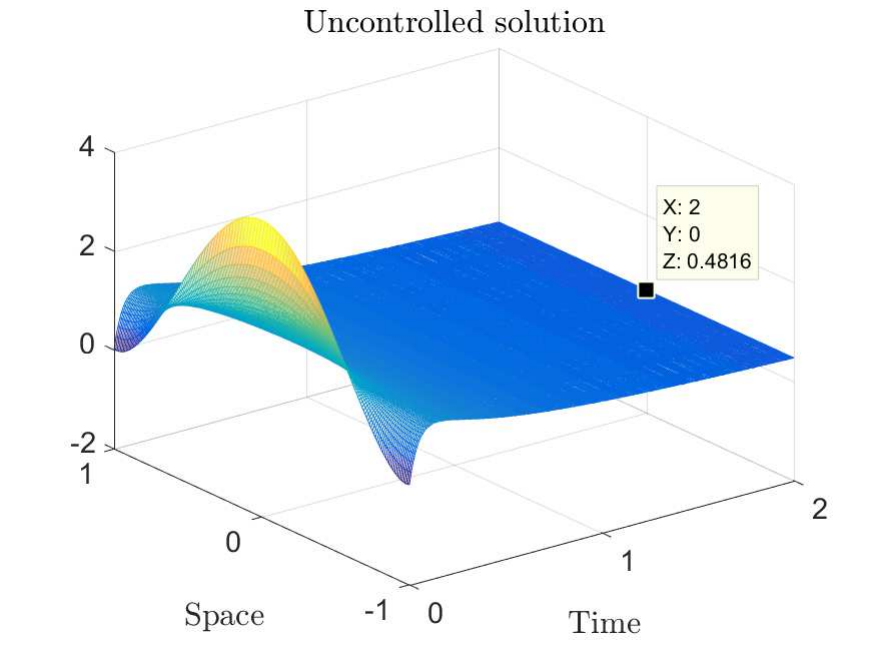}
\includegraphics[width=0.45\textwidth]{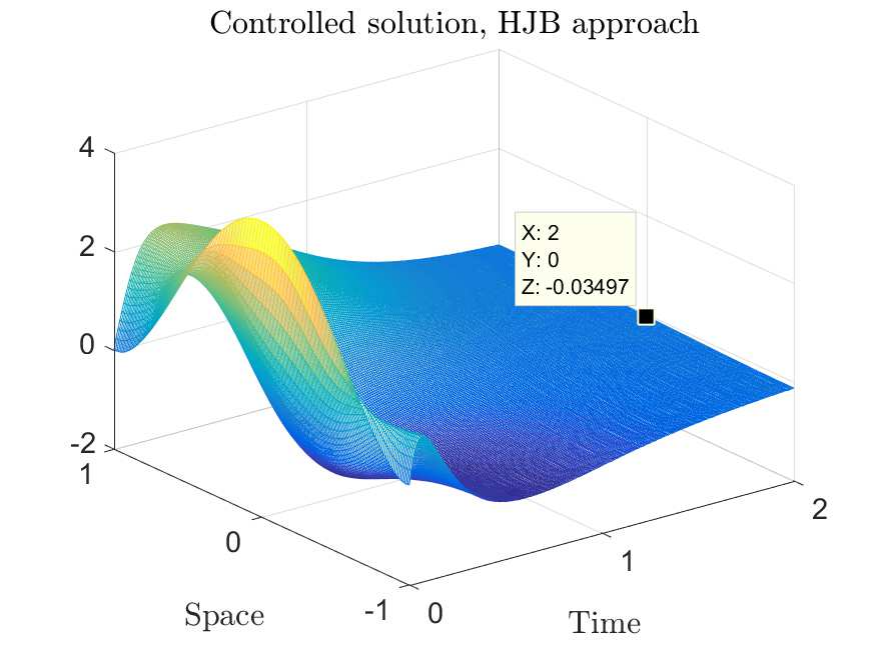}
\includegraphics[width=0.45\textwidth]{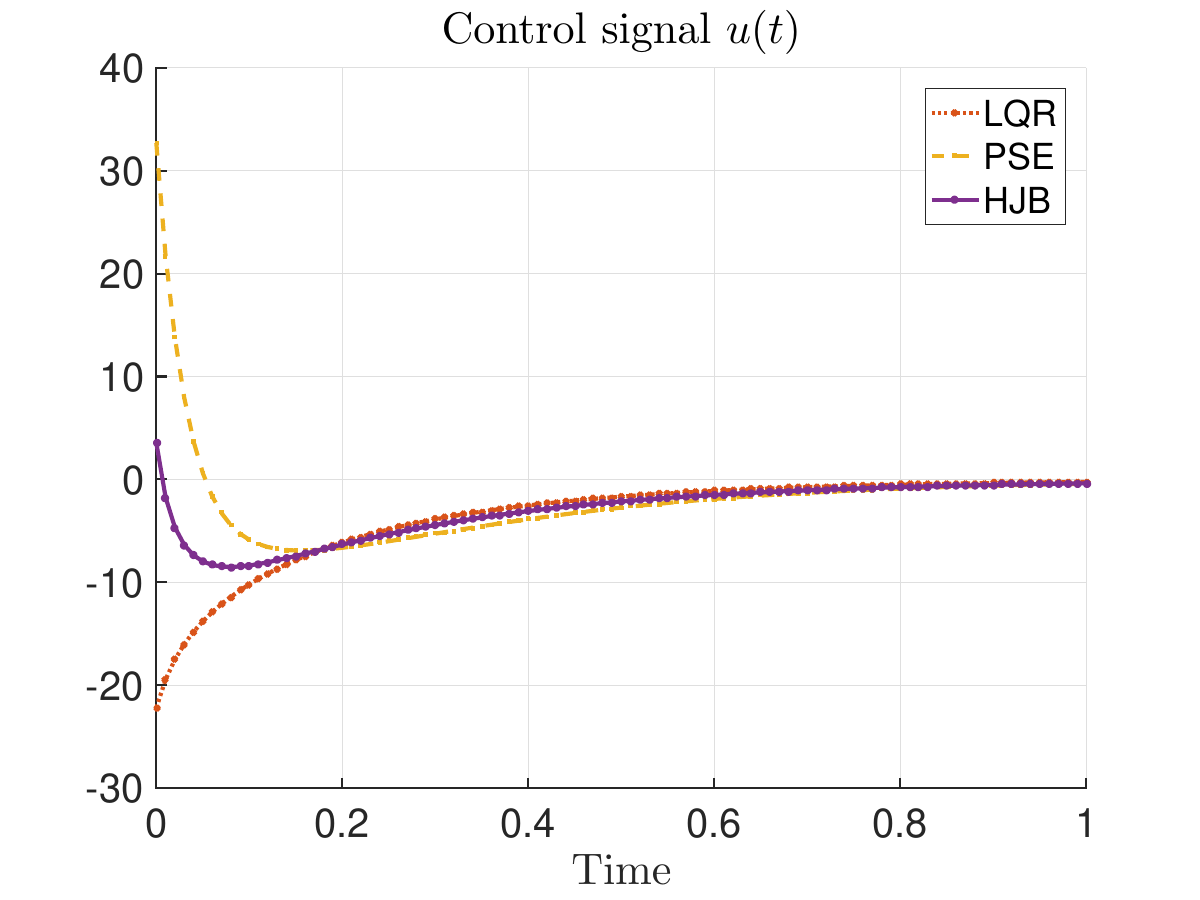}
\includegraphics[width=0.45\textwidth]{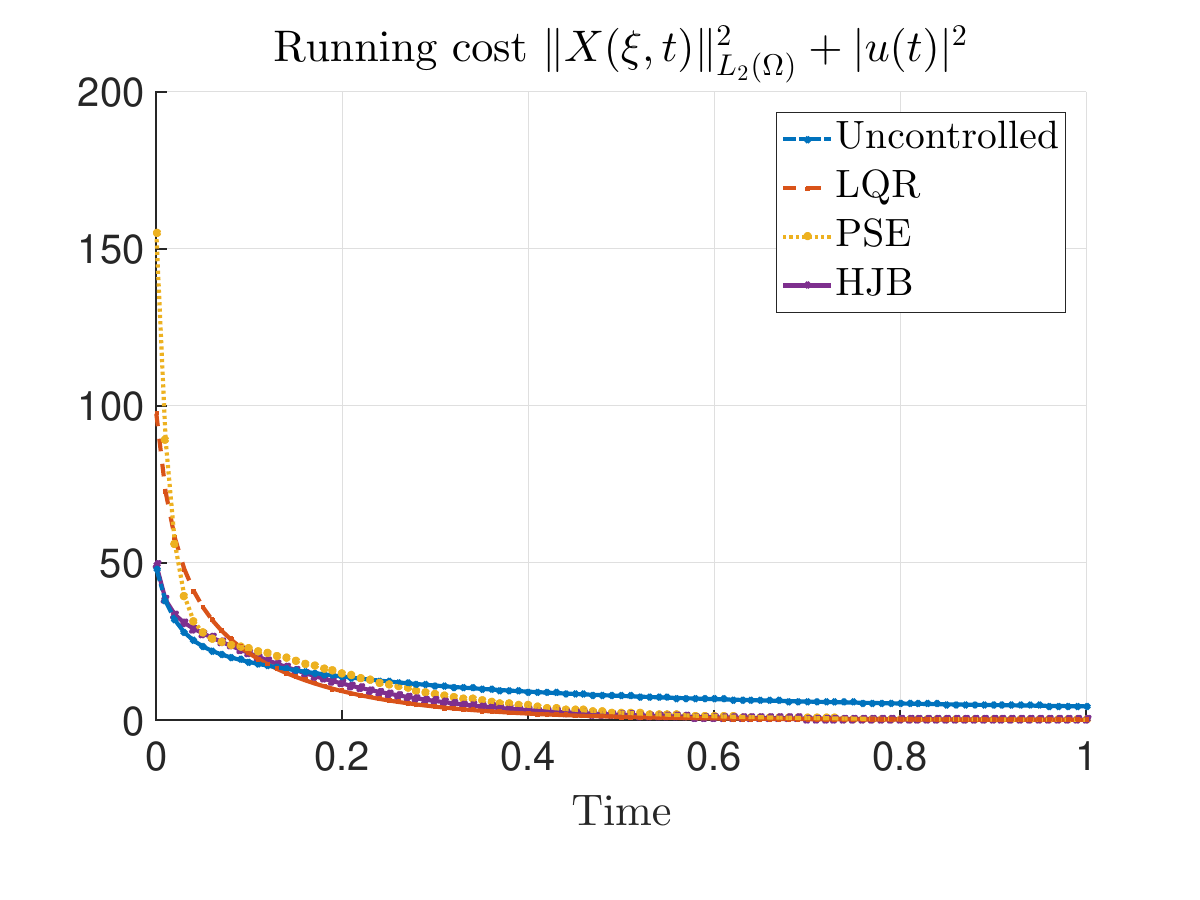}
\caption{A first preview of the stabilization of the semilinear parabolic equation \eqref{ex1}. Initial condition: $X_0(\xi)=4(\xi-1)^2(\xi+1)^2$. Dynamics are stable but slow. Total closed-loop costs $\cJ(u,X_0)$: {\bf i)} Uncontrolled: 13.45, {\bf ii)} LQR: 7.39, {\bf iii)} PSE: 9.43, {\bf iv)} HJB: 6.56 .}\label{fig:ex1}
\end{figure}

We now present a preview of the numerical results of the proposed approach. Further details of the numerical scheme will be developed in the forthcoming sections. The system dynamics in \eqref{ex1}, are approximated in 12 collocation points (14 with b.c.s'), and therefore our approximation scheme seeks for a solution of a 12-dimensional HJB equation, which allows the computation of online optimal feedback controllers. We compare our HJB-based controller (HJB) to the linear-quadratic controller (LQR) obtained by linearization of the system dynamics, and to an approximation method for the HJB equation based on power series expansion (PSE) \cite{G72,TBR10}. In Figure \ref{fig:ex1} we observe the basic features of the dynamics and the control schemes. The uncontrolled system dynamics (diffusion+dissipative source term) are stable, but stabilization is extremely slow. The control algorithms considerably reduce the transient phase. However, the control signals are different, and the HJB-based controller generates a feedback control with reduced overall cost \eqref{costex1}. Observe that at the beginning of the time horizon even the signs of the LQR-, PSE-, and HJB-based controls differ.

\section{Approximate iterative solution of HJB equations}
In this section, we construct a numerical scheme for the approximation of the HJB equation
\begin{equation}\label{hjbs3}
\underset{u\in U}{min}\{ DV(x)^t(f(x)+g u)+ \ell(x)+\gamma|u|^2\}=0\,,\quad V(0)=0\,,
\end{equation}
where $U\equiv \R$. We recall two additional features in this equation which render the application of classical approximation techniques difficult: the absence of a variational formulation, and  the minimization with respect to the control variable $u$, which makes the HJB equation fully nonlinear. The simplest numerical approach to these problems is the use of monotone, grid-based discretizations (finite differences, semi-Lagrangian), in conjunction with a fixed point iteration for the value function $V,$ which typically depends on the use of a discount factor. The so-called ``value iteration'' procedure was first presented by Bellman in \cite{B57}, and although it has become a standard solution method for low-dimensional HJB equations, it suffers from three major drawbacks. First, the grid-based character of the scheme makes it inapplicable for high-dimensional dynamics, as the total number of degrees of freedom scales exponentially with respect to the dimension of the dynamical system. This corresponds to the most classical statement of the so-called {\em curse of dimensionality}. Second, the contractive mapping includes a minimization procedure which needs to be solved for every grid point at every iteration. Third, the Lipschitz constant of the contractive mapping goes to 1 when the discretization parameter goes to 0, becoming extremely slow for fine-mesh solutions. In order to circumvent these limitations, we develop a numerical scheme combining an iteration on the control variable rather than the value function, together with a polynomial expansion for the value function to mitigate the computational burden associated to mesh-based schemes.
\subsection{Successive approximation of HJB equations}
In the following, we revisit the method presented in \cite{BST7,BST98}, which is referred as Successive Approximation Algorithm. We begin by defining the set of admissible controls.
\begin{definition}[Admissible control] We say that a feedback mapping $u:=u(x)$ is admissible on $\Omega\subset\R^d$, denoted as $u\in\cA(\Omega)$, if $u(x)\in \cC(\Omega)$, $u(0)=0$, and $\cJ(u(x(\cdot)),x_0)<\infty$ for all $x_0\in\Omega$.
\end{definition}
Starting from an admissible initial guess $u^0(x)$, the Successive Approximation Algorithm (Algorithm \ref{alg:sga1} below) generates the pair $(V^*,u^*)$ which solves equation \eqref{hjbs3}.
\begin{algorithm}[!ht]
\begin{algorithmic}
\STATE{Given $u^0(x)\in\cA(\Omega)$ and $tol>0$ \;}
\WHILE{$error>tol$}
\STATE{Solve
\begin{equation}
DV^i(x)^t(f(x)+g u^i)+ \ell(x)+\gamma|u^i|^2=0\,,\quad V^i(0)=0\,.
\end{equation}}
\STATE{Update
\[
u^{i+1}(x)=-\frac{1}{2\gamma}g^tDV^i(x)\,,
\]
$error=\|V^{i}-V^{i-1}\|$}
\ENDWHILE
\RETURN{$(V^*,u^*)$}
\caption{Successive Approximation Algorithm}\label{alg:sga1}
\end{algorithmic}
\end{algorithm}
Algorithm \ref{alg:sga1} corresponds to a Newton method for solving equation \eqref{hjbs3}, and in the linear-quadratic setting it is equivalent to the Newton-Kleinmann iteration for solving the Riccati equation. It can be also directly identified with the policy iteration algorithm for HJB equations (see \cite{AFK15} and references therein), although in this context the usual setting includes a discount factor which relaxes the admissibility assumption, as well as discrete-time dynamics. Consequently, it is applied to a Bellman equation with no continuous gradient. In both cases, the core ingredient of the algorithm is to generate a decreasing sequence of values $V^i$ by solving an associated sequence of linear problems. In our case this translates into solving, for a given $u(x)$ at each iteration, the Generalized Hamilton-Jacobi-Bellman (GHJB) equation
\begin{align}\label{ghjb}
\cG(DV;u)=&0\,,\quad V(0)=0\,,\\
\cG(p,u):=&p^t(f(x)+g u)+ \ell(x)+\gamma|u|^2\,.\notag
\end{align}
The following result from \cite{BST7} summarizes relevant properties of the GHJB equation.

\begin{proposition}\label{tghjb}
If $\Omega$ is a compact subset of $\R^d$, $f(x)$ is Lipschitz continuous on $\Omega$ and $f(0)=0$,  $l(x)\geq 0$ is strictly increasing in $\Omega$, $\gamma>0$, and $u\in\cA(\Omega)$, then:
\begin{enumerate}
\item There exists a unique $V(x)\in C^1(\Omega)$ satisfying \eqref{ghjb}.
\item $V(x)$ is a Lyapunov function of the controlled system.
\item $V(x)=\cJ(u,x)$, for all $x\in\Omega$.
\item The update $u^{+}(x):=-\frac{1}{2\gamma}g^tDV(x)$ satisfies $u^{+}\in\cA(\Omega)$.
\item If $V^+$ satisfies $\cG(DV^+;u^+)=0$, then $V^+\leq V$ for all $x\in\Omega$.
\end{enumerate}
\end{proposition}

\subsection{A continuation procedure} A critical aspect of the Successive Approximation Algorithm \ref{alg:sga1} is its initialization, which requires the existence of an admissible control $u^0(x)$ which in view of \eqref{costex1} means that it asymptotically stabilizes all the initial conditions in $\Omega$. For asymptotically stable dynamics, this is trivially satisfied by $u^0(x)=0$. For more general cases, the computation of stabilizing feedback controllers is a challenging task. A partial answer is to consider the stabilizing feedback associated to the linearized system dynamics. However, this feedback is only locally stabilizing, and therefore the identification of  a suitable domain $\Omega$ where this control law is admissible becomes relevant. For low dimensional dynamics, this has been studied in the context of Zubov's method in \cite{CGW01}. An alternative solution that we propose is to consider a discounted infinite horizon control problem

\[\underset{u(\cdot)\in\cU}{\min}\;\cJ(u(\cdot),x_0):=\int\limits_0^\infty e^{-\lambda t}\,(\ell(x(t))+\gamma|u(t)|^2)\, dt\,,\qquad\lambda>0\,,
\]
where the inclusion of the discount factor $\lambda$ relaxes the admissibility condition. Recently, in  \cite{GGKW16,PBND14}, the link between discounted optimal control and asymptotic stabilization has been discussed, and under certain conditions, the discounted control problem can generate optimal controls that are also admissible for the undiscounted problem. We recall that the associated HJB equation for the infinite horizon optimal control problem is given by
\begin{equation}\label{hjbd}
\lambda V(x)+\underset{u\in U}{min}\{ DV(x)^t(f(x)+g u)+ \ell(x)+\gamma|u|^2\}=0\,,\quad V(0)=0\,,
\end{equation}
and the associated GHJB reads
\begin{align}\label{ghjbd}
\cG_{\lambda}(V,DV;u)=&0\,,\quad V(0)=0\,,\\
\cG_{\lambda}(q,p,u):=&\lambda q+p^t(f(x)+g u)+ \ell(x)+\gamma|u|^2\,.\notag
\end{align}
We consequently modify the Successive Approximation Algorithm in order to embed it within a path-following iteration with respect to the discount factor:

\begin{algorithm}[H]
\begin{algorithmic}
\STATE{Given $\lambda>0$, $\epsilon>0$, and $\beta \in (0,1)$,}
\WHILE{$\lambda>\epsilon$}
\STATE{Solve for $(V,u)$
\begin{equation}\label{eq:alg2}
\lambda V(x)+\underset{u\in U}{min}\{ DV(x)^t(f(x)+g u)+ \ell(x)+\gamma|u|^2\}=0\,,
\end{equation}
with Algorithm \ref{alg:sga1} and initial guess $u^0$.}
\STATE{Update
\begin{align*}
u^{0}&=u\,,\\
\lambda&=\beta\lambda\,.
\end{align*}
}
\ENDWHILE
\RETURN{$(V^*,u^*)$}
\caption{A Discounted Path-Following Approximation Algorithm}\label{alg:sga2}
\end{algorithmic}
\end{algorithm}

For a sufficiently large $\lambda$, this algorithm can be initialized with $u_{\lambda}^0=0$. Continued reduction of the discount factor using hotstart every time when \eqref{eq:alg2} is called with a reduced $\lambda$-value, leads to an approximate solution of equation \eqref{hjbs3}.

\subsection{Spectral element approximation of the GHJB equation}
So far we have discussed the iterative aspects of a computational method for solving HJB equations.  We now address the numerical approximation of the GHJB equation.
\begin{equation}
\cG_{\lambda}(V,DV;u)=0\,,\quad V(0)=0\,.
\end{equation}
For this purpose, we consider an expansion $V_n(x)$ of the form
\[V_n(x)=\sum_{j=1}^n c_j\phi_j(x)\equiv\Phi_n\bc \,,\]
where $\Phi_n:=(\phi_1(x),\ldots,\phi_n(x))$, with $\phi_j\in C^{\infty}(\Omega,\R)$ belonging to a complete set of basis functions in $L^2(\Omega,\R)$, and $\bc=(c_1,\ldots,c_n)^t$. In particular, we shall often generate $\Phi_n$ from a multidimensional monomial basis as illustrated in Figure \ref{fig:pol2d}, which directly satisfies the boundary condition $V_n(0)=0$. The coefficients $c_j$ are obtained by imposing the Galerkin residual equation
\begin{equation}\label{ghjbg}
\langle \cG_{\lambda}(V_n,DV_n;u),\phi_i\rangle_{L^2(\Omega)}=0\,,\quad \forall \phi_i\in\Phi_n\,.
\end{equation}

\begin{remark}
The convergence of $V_n$ has been studied thoroughly in \cite{BST7}. It follows a power series argument, and requires conditions for uniform convergence of pointwise convergent series, in order to guarantee that $u_n:=-\frac12 \gamma^{-1}g^t D V_n(x)\in\cA(\Omega)$ for $n$ sufficiently large. In our particular case, we further assume that the dynamics $(f,g)$ are polynomial (as illustrated in Section \ref{pdesub}). Therefore, under the assumptions of Theorem 26 in \cite{BST7}, by choosing a multidimensional monomial basis  (of degree $\geq 2$) and an admissible control $u^0\in\cA(\Omega)$,  it can be established that, $\forall\epsilon>0$, $\exists K$ such that for $n>K$, $\|V-V_n\|_{L^2(\Omega)}<\epsilon$, and $u_n(x)\in\cA(\Omega)$.
\end{remark}

We now focus on the different terms involved in the approximation of the GHJB equation. Since this equation is meant to be solved within the iterative loop described in the previous section, we assume that $u(x)$ can be expressed in the form
\begin{equation}\label{u0}
u(x)=-\frac12 \gamma^{-1}g^t D V_n^0(x)\,,
\end{equation}
where $V^0(x)$ corresponds to the value function of the previous iteration, approximated with the expansion
\[V_n^0(x)= \sum_{j=1}^n c_j^0\phi_j(x).\]
Below we shall write $\bc^0$ for $(c^0_1,\dots,c^0_n)^t$. We proceed by expanding case by case the different terms of the Galerkin residual equation
\begin{equation}\label{ghjbge}
\langle \lambda V_n+DV_n^t(f(x)+g u)+ \ell(x)+\gamma|u|^2,\phi_i\rangle_{L^2(\Omega)}=0\,,\quad \forall \phi_i\in\Phi_n\,.
\end{equation}

\begin{enumerate}

\item[\bf 1)] \noindent$\pmb{\la\lambda V_n,\phi_i\ra_{L^2(\Omega)}}$: it is directly verifiable that
\[\la\lambda V_n,\phi_i\ra_{L^2(\Omega)}=\bM_{(i,\bullet)} \bc\,,\quad \bM\in\R^{n\times n}\,,\quad \bM_{(i,j)}=\lambda\la\phi_i,\phi_j\ra_{L^2(\Omega)}\,.\]

\item[\bf 2)] \noindent$\pmb{\langle  DV_n^tf,\phi_i\rangle_{L^2(\Omega)}}$: by inserting the expansion we obtain
\[DV_n^tf=\sum_{j=1}^n c_jD \phi_j^tf\,,\]
and therefore
\[\langle D V_n^tf,\phi_i\rangle_{L^2(\Omega)}=\bF_{(i,\bullet)} \bc\,\,,\quad \bF\in\R^{n\times n}\,,\quad\bF_{(i,j)}:=\langle D\phi_j^tf,\phi_i\rangle_{L^2(\Omega)}\,.\]

\item[\bf 3)] \noindent$\pmb{\la D V_n^t gu,\phi_i\ra_{L^2(\Omega)}}$: the relation \eqref{u0} leads to
\[
D V_n^t gu=D V_n^t\left(-\frac12 \gamma^{-1} gg^tD V_n^0\right)=-\frac12 \gamma^{-1} \sum_{j=1}^nc_jD\phi_j^t\left( gg^t\sum_{k=1}^nc_k^0D\phi_k\right)^t\,,
\]
such that
\begin{align*}
\la D V_n^t gu,\phi_i\ra_{L^2(\Omega)}&=\bG_{(i,\bullet)} \bc\,,\quad \bG\in\R^{n\times n}\,,\\
\bG_{(i,j)}&=-\frac12 \gamma^{-1}\sum_{k=1}^nc_k^0 \la g^tD\phi_kD\phi_j^t g,\phi_i\ra_{L^2(\Omega)}\,.
\end{align*}
\item[\bf 4)] \noindent$\pmb{\langle l(x),\phi_i\rangle_{L^2(\Omega)}}$: we further assume that
\[\langle l(x),\phi_i\rangle_{L^2(\Omega)}=\langle x^t Q x,\phi_i\rangle_{L^2(\Omega)}\,,\qquad Q> 0\in \R^{d\times d}\,.\]

\item[\bf 5)] \noindent$\pmb{\langle\gamma|u|^2,\phi_i\rangle_{L^2(\Omega)}}$ : note that
\[
\gamma|u|^2=\frac14 \gamma^{-1}(g^tD V_n^0)^ 2=\frac14 \gamma^{-1}\left(\sum_{j=1}^n c_j^0 g^tD\phi_j\right)^ 2\,,\]
leading to
\[\langle\gamma|u|^2,\phi_i\rangle_{L^2(\Omega)}=(\bc^0)^t\bU_{(i,\bullet)}\bc^0\,,\]
 $\bU\in\R^{n\times n\times n}$ is given by
\[\bU_{(i,j,k)}=\langle(g^tD\phi_{j})(g^tD\phi_{k}),\phi_{i}\rangle_{L^2(\Omega)}\,.\]
\end{enumerate}
After discretization, the GHJB \eqref{ghjbg} reduces to a parameter-dependent linear system for $\bc$
\[\left(\bM+\bF+\bG(\bc^0)\right)\bc=\bb(\bU,\bc^0)\,,\]
where $\bb$ is given by the expansion of $l(x)+\gamma|u|^2$ ( terms {\bf 4)} and {\bf 5)} in the list above).

\section{Computation of integrals via separable expansions}\label{sec:comp}

\begin{figure}[!h]
\centering
\includegraphics[width=0.8\textwidth]{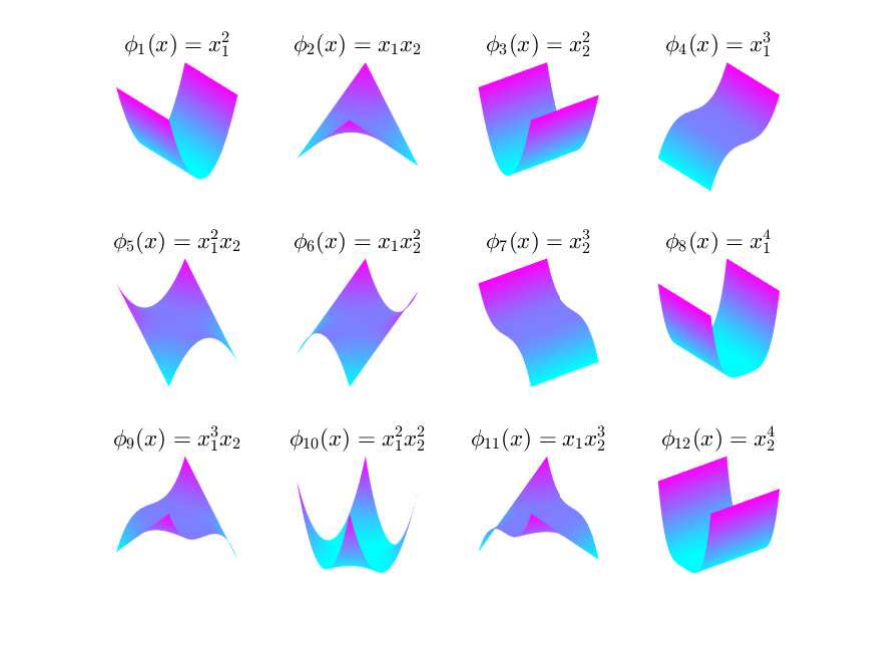}
\caption{Two dimensional monomial basis. The first three basis functions correspond to the terms of the Riccati ansatz for the linear-quadratic control problems, where the value function is known to be a quadratic form $x^t\Pi x$. Adding terms of higher order allows a more accurate solution for nonlinear control problems. We construct the high-order terms by limiting the degree of the monomials.}\label{fig:pol2d}
\end{figure}

 Under Assumption \ref{as:1} concerning the separability of  the free dynamics $f$, and with the construction of a separable set of basis functions by taking the tensor product of one-dimensional basis functions as shown in Figure \ref{fig:pol2d}, the calculation of the $d$-dimensional inner products of the Galerkin residual equation of the previous section is reduced to the product of one-dimensional integrals. In the following, we provide further details of this procedure.

\subsection{Generation of a multi-dimensional basis} The multi-dimensional basis functions \\$\Phi_n:=(\phi_1(x),\ldots,\phi_n(x))$ for the expansion of $V_n$ are generated as follows. We start by choosing a polynomial degree $M\in\mathbb{N}$, and a one-dimensional  polynomial basis $\varphi_M:\R\to \R^M$. For the sake of simplicity, we consider the monomial basis $\varphi_M=(1,x,\ldots,x^M)^T$, but the same ideas apply for other basis, such as orthogonal polynomials. The multidimensional basis is generated as a subset of the $d$-dimensional tensor product of one-dimensional basis, such that
\[\Phi_n\equiv\left\{\phi\in\bigotimes\limits_{i=1}^d\varphi_M(x_i)\,,\;\text{and } deg(\phi)\leq M\right\}\,\]
i.e., we construct a full multidimensional tensorial basis and then we remove elements according to the approximation degree $M$. The elimination step is fundamental and is twofold. If no elimination is performed, the cardinality of $\Phi_n$ would be $M^d$, and again one would face the \textsl{curse of dimensionality} that also affects grid-based schemes. By reducing the set to multdimensional monomials of degree at most $M$, the cardinality $n$ of the set $\Phi_n$ is given by
\begin{equation}\label{totalcount} n=\sum\limits_{m=1}^M\left(\begin{array}{c}d+m-1\\m\end{array}\right)\,,
\end{equation}
 which replaces the exponential dependence on $d$ by  a combinatorial one. This formula is evaluated in Table \ref{tab:count} for different values of interest for $M$ and $d$. By considering globally defined polynomial basis functions, the dependence on the dimension is replaced by the combinatorial expression \eqref{totalcount}. The dimensional reduction of the basis is particularly significant for low order polynomial approximation (up to degree 6).
 A second justification for the way in which we generate the basis set has a control-theoretical inspiration. A well-known result in optimal feedback control is that if the dynamics are linear, and the running cost is quadratic, the value function associated to the infinite horizon control problem (in the unconstrained case and other technical assumptions) is a quadratic form, i.e. is of the form $V(x)=x^t\Pi x$, which fits precisely the elements generated for $\Phi_n$ with a monomial basis when $M=2$ and linear elements are eliminated. Therefore, our basis can be interpreted as a controlled increment, accounting for the nonlinear dynamics, of the basis required to recover the solution of the control problem associated to the linearized dynamics around the equilibrium point.

\begin{table}[h]
\centering
\begin{tabular}{|c|cccc|cccc|}
\hline
& \multicolumn{4}{c}{Full monomial basis} &\multicolumn{4}{c|}{Even-degree monomials}\\
\hline
$d$\textbackslash$M$& 2 & 4 & 6& 8& 2 & 4 & 6& 8\\
\hline
\rule{0pt}{4ex}
6&27 & 209 &923&3002&21 & 147 &609&1896\\
8& 44& 494&3002&12869&36 & 366 &2082&8517\\
10& 65& 1000&8007&43757&55 & 770 &5775&30085\\
12& 90& 1819&18563&125969&78 & 1443 &13819&89401\\
14&  119 & 3059 & 38759& 319769 &105 &2485 & 29617 &233107\\
\hline
\end{tabular}
\vskip 2mm
\caption{Number of elements $n$ in the basis, as a function of the dimension $d$ and the total polynomial degree $M$. The global polynomial approximation partially circumvents the curse of dimensionality, as the dimension of the basis no longer depends exponentially on the dimension, but rather combinatorially.}\label{tab:count}
\end{table}

\begin{remark}\label{oddeven}
Theorem 7.1 in \cite{BST98} states parity conditions to reduce the polynomial basis $\Phi_n$. Under the assumptions $l(x)=x^t Q x$, and $g\in\R^d$,  if
\begin{itemize}
\item[i)] $\Omega$ is a symmetric rectangle around the origin, i.e., $\Omega=[-l_1,l_1]\times\ldots\times[-l_d,l_d]\,,$
\item[ii)] the free dynamics are odd-symmetric on $\Omega$, i.e. $f(-x)=-f(x)$, for all $x\in\Omega\,,$
\end{itemize}
then $V_n(x)$ is an even-symmetric function, i.e., $V_n(-x)=V_n(x)$, and therefore odd-degree monomials are excluded from the basis. A direct corollary is that in the linear quadratic case, where the linear dynamics are trivially odd-symmetric, $V(x)$ is a quadratic form.
\end{remark}

\noindent Finally, for the calculation presented in the following, it is important to note that due to the construction procedure, the basis elements directly admit a separable representation
\begin{equation}\label{eq:4.0}
\phi_i(x)=\prod_{j=1}^d \phi_i^j(x_j)=\prod_{j=1}^dx_j^{\nu_j}\,,\;\;\text{with } \sum_j\nu_j\leq M\,,
\end{equation}
where each component $\phi_i^j(x)\in\varphi_M$.

\subsection{High-dimensional integration}
We begin by recalling that
\begin{equation}\label{eq:4.1}
f_i(x)=\sum_{j=1}^{n_f}\prod_{k=1}^{d}\cF_{(i,j,k)}(x_k)\,,
\end{equation}
where $\cF(x):\R^d\rightarrow\R^{d\times n_f\times d}$ is a tensor-valued function, and that $g\in\R^d$.

As in the previous section, we proceed term by term, to obtain the summands in \eqref{ghjbge}. The integration is carried over the hyperrectangle $\Omega=\Omega_1\times\ldots\times\Omega_d$.

\begin{enumerate}

\item[\bf 1)] \noindent$\pmb{\la\lambda V_n,\phi_i\ra_{L^2(\Omega)}}$: this term is directly assembled from the calculation of
\[\la\phi_i,\phi_j\ra_{L^2(\Omega)}=\prod_{k=1}^d\int_{\Omega_k} \phi_i^k(x_k)\phi_j^k(x_k)\,dx_k\]

\item[\bf 2)] \noindent$\pmb{\langle  DV_n^tf,\phi_i\rangle_{L^2(\Omega)}}$: This term involves the calculation of
\[\langle D\phi_j^tf,\phi_i\rangle_{L^2(\Omega)}=\sum_{p=1}^d\langle f_p\px{p}\phi_j,\phi_i\rangle_{L^2(\Omega)}\,.\]
which is expanded by using the separable structure of the free dynamics
\[\langle f_p\px{p}\phi_j,\phi_i\rangle _{L^2(\Omega)}=\sum\limits_{l=1}^{n_f}\langle \left(\prod_{m=1}^d\cF(p,l,m)\right)\px{p}\phi_j,\phi_i\rangle_{L^2(\Omega)}\,,\]
where
\begin{align*}
&\langle\left(\prod_{m=1}^d\cF(p,l,m)\right)\px{p}\phi_j,\phi_i\rangle_{L^2(\Omega)}\\ 
=&\left(\prod_{\begin{subarray}{l}m=1\\m\neq p\end{subarray}}^d\int\limits_{\Omega_m}\cF(p,l,m)\phi_i^m\phi_j^m(x_m)\,dx_m\right)\left(\int\limits_{\Omega_p}
\cF(p,l,p)\phi_i^p\px{p}\phi_j^p(x_p)\,dx_p\right)
\end{align*}

\item[\bf 3)] \noindent$\pmb{\la D V_n^t gu,\phi_i\ra_{L^2(\Omega)}}$: In this case, we need to work on the expression
\[\langle g^t D\phi_kD\phi_j^t g,\phi_i\rangle_{L^2(\Omega)}=\sum\limits_{l,m=1}^dg_lg_m\langle \px{l}\phi_k\px{m}\phi_j,\phi_i\rangle_{L^2(\Omega)},\]
which is obtained directly from the computations for $\pmb{\langle\gamma|u|^2,\phi_i\rangle_{L^2(\Omega)}}$ in {\bf{5)}} below.

\item[\bf 4)] \noindent$\pmb{\langle l(x),\phi_i\rangle_{L^2(\Omega)}}$: \[\langle l(x),\phi_i\rangle_{L^2(\Omega)}=\langle x^t Q x,\phi_i\rangle_{L^2(\Omega)}=\sum\limits_{j,k=1}^dQ_{(j,k)}\langle x_jx_k,\phi_i\rangle_{L^2(\Omega)}\,,\]
where with a similar argument as in the previous term we expand
\[\langle x_jx_k,\phi_i\rangle_{L^2(\Omega)}=\left(\prod_{\begin{subarray}{l}p=1\\p\neq j\\p\neq k \end{subarray}}^d\int\limits_{\Omega_p}\!\!\phi_i^p(x_p)\,dx_p\right)\!\!\left(\int\limits_{\Omega_j}\!\!\phi_i^j(x_j)x_j\,dx_j\right)\!\!
\left(\int\limits_{\Omega_k}\!\!\phi_i^k(x_k)x_k\,dx_k\right)\,.\]

\item[\bf 5)] \noindent$\pmb{\langle\gamma|u|^2,\phi_i\rangle_{L^2(\Omega)}}$ : This term requires the computation of the inner product
\[\langle(g^ tD\phi_{j})(g^tD\phi_{k}),\phi_{i}\rangle_{L^2(\Omega)}=g^t \tilde\bU_{(\cI,\bullet)} g\,,\quad \cI=(i,j,k)\,,\]
with $\tilde\bU\in \R^{n\times n\times n\times d\times d}$ given by
\[\tilde\bU_{(\cI,l,m)}:=\langle\px{l}\phi_{j}\px{m}\phi_{k},\phi_{i}\rangle_{L^2(\Omega)}\,.\]
By using the separable representation of the basis functions
\[\px{l}\phi_{j}=\left(\prod_{\begin{subarray}{l}p=1\\ p\neq l\end{subarray}}^d\phi_j^p(x_p)\px{l}\right)\phi_j^l(x_l)\]
we expand the inner product
\begin{eqnarray}\label{eq:bU}
\tilde\bU_{(\cI,l,m)}&=&\left(\prod_{\begin{subarray}{l}p=1\\ p\neq l\\p\neq m\end{subarray}}^d\int\limits_{\Omega_p}\!\!\phi_i^p\phi_j^p\phi_k^p(x_p)\,dx_p\right)\!\!\left(\int\limits_{\Omega_l}\!\!\phi_i^l
\phi_k^l\px{l}\phi_j^l(x_l)\,dx_l\right)\!\!\ldots\\
&&\ldots\left(\int\limits_{\Omega_m}\!\!\phi_i^m\phi_j^m\px{m}\phi_k^m(x_m)\,dx_m\right)\,.\nonumber
\end{eqnarray}
\end{enumerate}

\paragraph{\bf Initialization} The first iteration, with a stabilizing initial guess $u^0$, requires special attention. If it is obtained via a Riccati-type argument, then initialization follows directly from \eqref{u0}. Otherwise we shall relax this requirement, and only assume that the initial stabilizing controller is given in separable form
\[u^0(x)=\sum_{j=1}^{n_u}\prod_{k=1}^{d}\cU^0_{(j,k)}(x_k)\,,\]

In this case, we must recompute the term:
\begin{itemize}
\item ${\pmb{\langle \gamma|u^0|^2,\phi_i\rangle_{L^2(\Omega)}}}$
\begin{align*}
\langle\gamma|u^0|^2,\phi_i\rangle_{L^2(\Omega)}&=\gamma\langle(\sum_{j=1}^{n_u}\prod_{k=1}^{d}\cU^0_{(j,k)}(x_k))^2,\phi_i\rangle\,\\
&=\gamma\sum\limits_{j,l=1}^{n_u}\langle\left(\prod_{k=1}^d\cU^0(j,k)\right)\left(\prod_{k=1}^d\cU^0(l,k)\right),\phi_i\rangle\,,\\
&=\gamma\sum\limits_{j,l=1}^{n_u}\prod_{k=1}^d\int\limits_{\Omega_k}\cU^0(j,k)\cU^0(l,k)\phi_i^k (x_k)\,dx_k\,.
\end{align*}

As for the term $\langle DV_n^t gu^0 ,\phi_i\rangle_{L^2(\Omega)}$, which needs to be computed differently in the first iteration, we can proceed in the same way as for  $\langle  D V_n^t f,\phi_i\rangle_{L^2(\Omega)}$, since both $gu^0$ and $f$ have the same separable structure, it just takes to assign $f_i=g_iu^0$.
\end{itemize}

\subsection{Computational complexity and implementation}
Among the expressions developed in the previous subsection, the overall computational burden is governed by the approximation of
\[\langle\gamma|u|^2,\phi_i\rangle_{L^2(\Omega)}\,,\]
which requires the assembly of the 5-dimensional tensor $\tilde\bU\in\R^{n\times n\times n\times d\times d}$. Each entry of this tensor is a $d$-dimensional inner product, which under the aforementioned separability assumptions is computed as the product of $d$, one-dimensional integrals. Thus, the total amount of one-dimensional integrals required for the proposed implementation is $O(n^3d^3)$.
\noindent A positive aspect of our approach is that the assembly of tensors like $\tilde\bU$ falls within the category of \textsl{embarrassingly parallelizable} computations, so the CPU time scales down almost directly with respect to the number of available cores.  Furthermore, $\tilde\bU$ can be entirely computed in an offline phase, before entering the iterative loops in Algorithms \ref{alg:sga1} and \ref{alg:sga2}. However, for values of interest of $n$ and $d$, such as $d>10$ and $n=4$, Table \ref{tab:count} indicates that $n^3d^3$ is indeed a very large number. A rough estimate of the CPU time required for the assembly of $\tilde\bU$ is given by
\[\text{CPU}(\tilde\bU)=\frac{t_{1d}\times n^3\times d^3}{\#cores}\,,\]
where $t_{1d}$ corresponds to the time required for the computation of a one-dimensional integral. Therefore, it is fundamental for an efficient implementation to reduce $t_{1d}$ to a bare minimum.
From closer inspection of the expression \eqref{eq:bU}, we observe that all the terms can be identified as elements of the tensors $\mathcal{M},\mathcal{K}\in\R^{M\times M\times M}$
\[
\mathcal{M}_{(i,j,k)}:=\int\limits_{x_l}^{x_u}\varphi_i(x)\varphi_j(x)\varphi_k(x)\,dx\,,\quad \mathcal{K}_{(i,j,k)}:=\int\limits_{x_l}^{x_u}\varphi_i(x)\varphi_j(x)\partial_x\varphi_k(x)\,dx\,.\]
Both $\mathcal{M}$ and $\mathcal{K}$ can be computed exactly with a Computer Algebra System, or approximated under suitable quadrature rules. We follow this latter approach, implementing an 8-point Gauss-Legendre quadrature rule. After having computed $\mathcal{M}$ and $\mathcal{K}$, the assembly of \eqref{eq:bU} reduces to $d$ calls to properly indexed elements of these tensors. This approach requires a careful bookkeeping of the separable components of each multidimensional basis function $\phi_i$. In this way, an entry of $\tilde\bU$ takes of the order of $10^{-7}$ seconds and the overall CPU time is kept within hours for problems of dimension up to 12.

\section{Computational implementation and numerical tests}

\subsection{Convergence of the polynomial approximation}
\begin{figure}[!h]
\centering
\includegraphics[width=.48\textwidth]{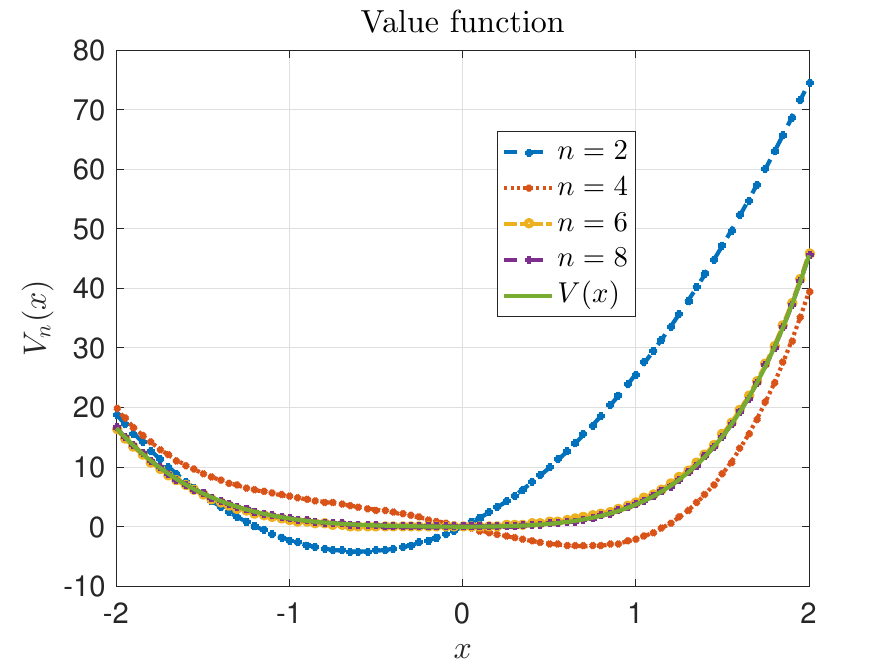}
\includegraphics[width=.48\textwidth]{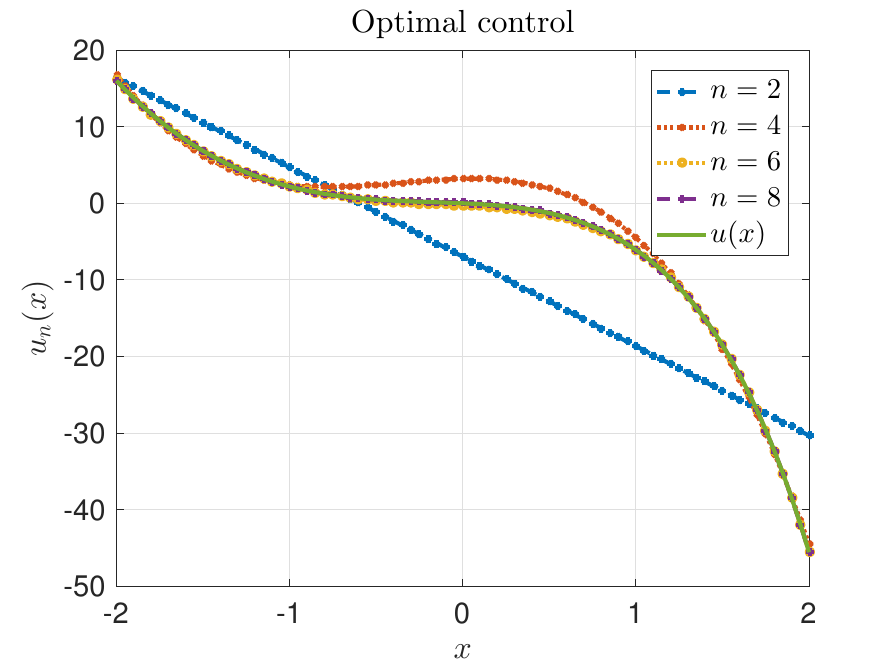}
\caption{1D polynomial approximation of the infinite horizon control problem with nonquadratic running cost. Approximation with monomial basis. The number $n$ denotes the total number of basis functions.}\label{figconv1d}
\end{figure}

We assess the convergence of the polynomial approximation in a 1D test, with
\[f(x)=0\,,\quad g=1\,,\quad l(x)=\frac{1}{4R}{\left(x^2\, \mathrm{e}^{x} + 2\, x\, \mathrm{e}^{x} + 4\, x^3\right)}^2\,,\quad\Omega=(-1,1)\,,\]
such that the exact solution of equation \eqref{hjb2} is given by
\[ V(x)=x^4+x^2e^x\,.\]
We implement the path-following version (Algorithm 2), starting with $u^0=0$, $\lambda=1$ and a threshold value $\epsilon=1\times 10^{-6}$, a parameter $\beta=0.5$, and an internal tolerance $tol=10^{-8}$. The relative error for Table \ref{tabtestconv} is defined as
\[\text{error}:=\frac{\|V_n(x)-V(x)\|_{L^2(\Omega)}}{\|V(x)\|_{L^2(\Omega)}}\]
 and number of iterations for different polynomial degree approximations are shown in Table \ref{tabtestconv} and Figure \ref{figconv1d}.

\begin{table}[h]
\centering
  \begin{tabular}{ccccc}
  \hline\\
  & \multicolumn{2}{c}{Monomial basis} & \multicolumn{2}{c}{ Legendre basis}\\
 \cmidrule(lr){2-3}\cmidrule(lr){4-5}\\
  $n$(degree) &error& iterations &error& iterations\\
\cmidrule(lr){1-1}\cmidrule(lr){2-2}\cmidrule(lr){3-3}\cmidrule(lr){4-4}\cmidrule(lr){5-5}\\
2 & 1.1539 & 53 &1.4127 & 52\\
4 & 0.2541 & 49 &0.3643 & 58\\
6 & 0.015 & 52 &0.0206 & 52\\
8 & 5.01$\times 10^{-4}$ & 55 &6.41$\times 10^{-4}$ & 53\\
10 & 8.33$\times 10^{-6}$ &55 &1.072$\times 10^{-5}$ & 55\\
    \hline
    \\
  \end{tabular}
  \caption{1D polynomial approximation of the infinite horizon control problem with nonquadratic running cost. The number $n$ denotes the total number of basis functions.}\label{tabtestconv}
\end{table}

\subsection{Optimal feedback control of semilinear parabolic equations}
Similarly as in Section \ref{pdesub}, we consider the following optimal control problem
\begin{equation}
\underset{u(\cdot)\in L^2([0;+\infty))}{\min} \cJ(u(\cdot),X_0):=\int_0^\infty\|X(\xi,t)\|_{L^2(\cI)}^2+\gamma u(t)^2\,dt\,,
\end{equation}
subject to the semilinear dynamics
\begin{align*}
\partial_{t}X(\xi,t) &=\mathcal{L}(X,X_{\xi\xi}) +\cN(X,\partial_\xi
X)+\chi_{\omega}(\xi)u(t)\,,\qquad\text{in}\;\cI\times\R^+\,,\\
%X_\xi(\xi_l,t)&=X_\xi(\xi_r,t)=0\,,\quad t\in \R^+\,,\\
X(\xi,0)&=X_0(\xi)\,,\quad \xi\in \cI\,,
\end{align*}
where the linear operator $\mathcal{L}$ is of the form $\mathcal{L}:=\sigma \partial_{\xi\xi}X(\xi,t)+rX(\xi,t)$ with $r\in\R$, and $\cN$ is a nonlinear operator such that $\cN(0,0)=0$. The scalar control acts through the indicator function $\chi_{\omega}(\xi)$, with $\omega\subset\cI$ The system is closed under suitable boundary conditions.
We choose $\cI=(-1,1)$, $\omega=(-0.5,-0.2)$, $\sigma=0.2$, and $\gamma=0.1$. The nonlinearity covers both advective, Burgers'-type, and polynomial source terms. In order to generate a low-dimensional state space representation of the dynamics, we resort to a pseudospectral collocation method with Chebyshev polynomials as in \cite{os06} (for further details we also refer to \cite[p. 107]{qvbook}. By considering $d$ collocation points $\xi_i= -cos(\pi i/d)\,,i=1,\ldots,d$, the continuous state $X(\xi,t)$ is discretized into $X(t)=(X_1(t),\ldots,X_d(t))^t\in\R^d$, where $X_i(t)=X(\xi_i,t)$. The semilinar PDE dynamics are thus approximated by the $d-$dimensional nonlinear ODE system
\begin{equation}
\dot X(t)=AX(t)+N(X(t))+Bu(t)\,,\label{eq:nlode}
\end{equation}
where the operators $(A,N,B)$ correspond to the finite-dimensional realization of  $(\mathcal{L},\cN,\chi_{\omega}(\xi))$ through the Chebyshev pseudospectral method. Therefore, the number of collocation points governs the dimension of the resulting nonlinear ODE system \eqref{eq:nlode}, and consequently determines the dimension of the domain $\Omega$ where the associated HJB equation is solved. In the following, Tests 1-3 are computed in 14 collocation points, which after including boundary conditions lead to a 12 dimensional domain $\Omega$ for the HJB equation.  Test 4 is solved in 14 dimensions. The high-dimensional solver was implemented in MATLAB, parallelizing the tensors assembly, and tests were run on a muti-core architecture 8x Intel Xeon E7-4870 with 2,4Ghz, 1 TB of RAM. The MATLAB pseudoparallelization distributes the tasks among 20 workers. Representative performance details are shown in Table \ref{tab:cpu}. The assembly of high-dimensional tensor that enter the iterative algorithm accounts for over 80\% of the total CPU time. This percentage increases when Algorithm \ref{alg:sga1} is implemented for asymptotically stable dynamics, as it requires a much lower number of iterations. Note that much of the work done during the assembly phase is independent of the dynamics (see for instance \eqref{eq:bU}), and therefore  can be re-used in latter problems, mitigating the overall computational burden.
\begin{table}[h]
\centering
\setlength{\tabcolsep}{1mm}
  \begin{tabular}{cccc}
  \hline\\
   Test &Dimension&CPU-assembly&CPU-iterative (\#)\\
   \cmidrule(lr){1-1}\cmidrule(lr){2-2}\cmidrule(lr){3-3}\cmidrule(lr){4-4}\\
   1 & 10 & 2.061$\times 10^{3}$[s] & 4.221$\times 10^{2}$[s](32)\\
   1 & 12 & 1.945$\times 10^{4}$[s] & 3.377$\times 10^{3}$[s](32)\\
   4 & 14 & 1.557$\times 10^{5}$[s] & 3.102$\times 10^{4}$[s](37)\\
    \hline
    \\
  \end{tabular}
  \caption{CPU times for different tests and dimensions. CPU-assembly corresponds to the amount of time spent in offline assembly of the different terms of the Galerkin residual equation \eqref{ghjbg}. CPU-iterative refers to the amount of time spent inside Algorithm \ref{alg:sga2}.}\label{tab:cpu}
\end{table}

We now turn to the specification of parameters for the solution of the HJB equation. We set $\Omega=(-2,2)^d$, and consider a monomial basis up to order 4 as described in Section \ref{sec:comp}. Depending on the dynamics of every example, we will neglect odd-degree basis functions as in Remark \ref{oddeven}. All the integrals are approximated with an 8 point Gauss-Legendre quadrature rule.
Whenever system dynamics are stable at the origin, the value function is obtained from the undiscounted Algorithm \ref{alg:sga1}, initialized with $u^0=0$. When the dynamics are unstable over $\Omega$, we implement Algorithm \ref{alg:sga2}, with $\lambda=1$, $\epsilon=10^{-6}$, and $\beta=0.9$. The initializing controller is given by the solution of the associated linear-quadratic optimal feedback, as described below. For both implementations, the tolerance of the algorithm is set to $tol=10^{-8}$.
In the following tests, we compare the HJB-based feedback control with respect to the uncontrolled dynamics ($u=0$), the linear-quadratic optimal feedback (LQR), and the power series expansion type of controller (PSE). We briefly describe these controllers.
The well-known LQR feedback controller corresponds to the HJB synthesis applied over the linearized system around the origin
\begin{equation}
\dot X(t)=AX(t)+Bu(t)\,,\label{eq:lode}
\end{equation}
and results in the optimal feedback control law given by
\[u^*=-\gamma^{-1}B^t\Pi X\,,\]
where $\Pi\in R^{d\times d}$  is the unique self-adjoint, positive-definite solution of  the algebraic  Riccati equation
\[
A^t\Pi+ \Pi A - \Pi B\gamma^{-1}B^t\Pi+Q=0\,,
\]
and $X^tQX$ corresponds to the finite-dimensional approximation of $\|X(\xi,t)\|_{L^2(\cI)}$.
Once this controller has been computed, the high-order PSE feedback is obtained as

\[u^*=-\gamma^{-1}B^t(\Pi X-(A^t-\Pi B\gamma^{-1}B^t)^{-1}\Pi N_l(X))\,,\]
where $N_l(X)$ corresponds to the lowest order term of the nonlinearity $N(X)$.Variations of such feedback laws have been discussed in previous publications, see eg. \cite{BK91} and references therein. For the Burgers' equation it was observed numerically  in  \cite{TBR10}  that this suboptimal nonlinear controller leads to an increased closed-loop stability region with respect to the LQR feedback applied for the linearized dynamics.

\subsection*{Test 1: Viscous Burgers'-like equation} In this first test we address nonlinear optimal stabilization of advective-reactive phenomena, by considering a 1D Burgers'-like model with $(\xi,t)\in\cI\times\R^+$ given by
\begin{align*}
\partial_{t}X(\xi,t) &=\sigma \partial_{\xi\xi}X(x,t) +X(\xi,t)\partial_{\xi}X(\xi,t)+1.5X(\xi,t)e^{-0.1X(\xi,t)}+\chi_{\omega}(\xi)a(t)\,,\\
X(\xi_l,t)&=X(\xi_r,t)=0\,,\quad t\in \R^+,\\
X(\xi,0)&=-\text{sign}(\xi)\,,\quad \xi\in \cI\,.
\end{align*}
\begin{figure}[!h]
\centering
\includegraphics[width=0.45\textwidth]{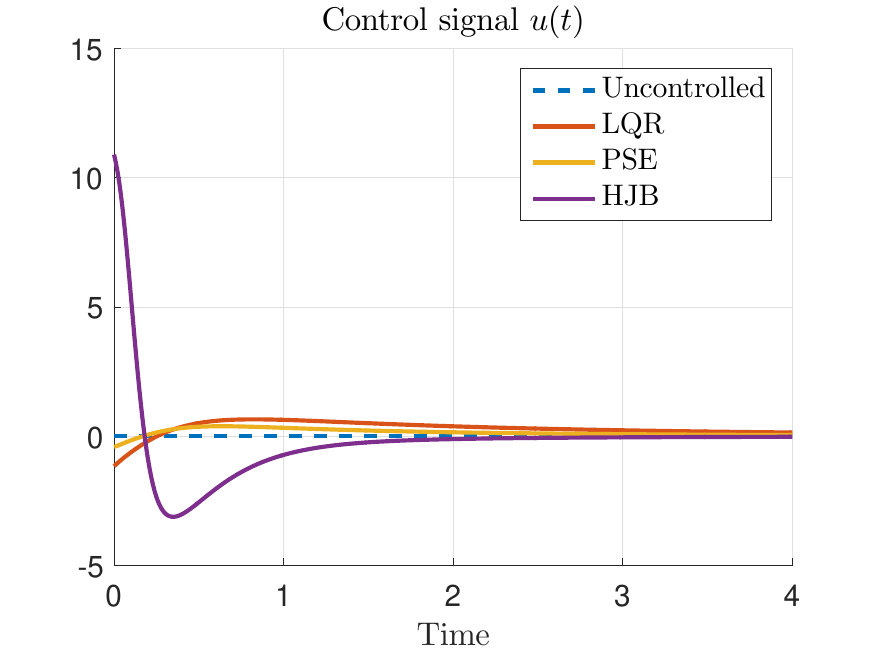}
\includegraphics[width=0.45\textwidth]{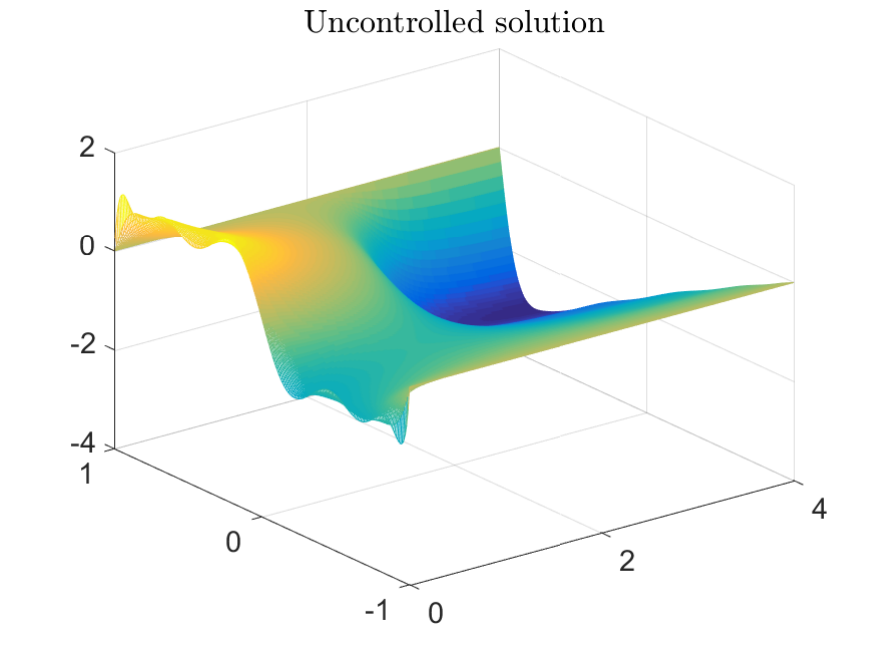}
\includegraphics[width=0.45\textwidth]{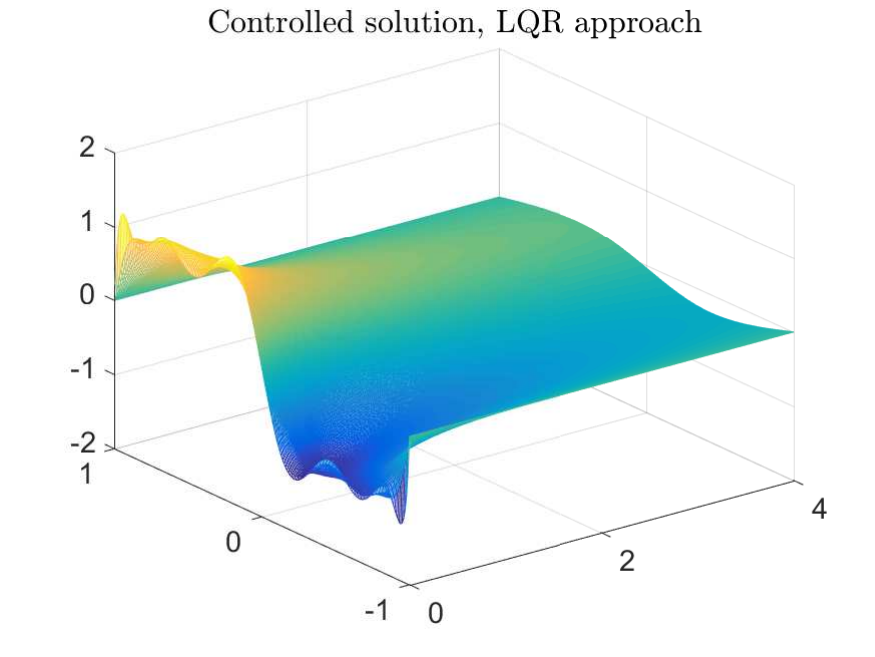}
\includegraphics[width=0.45\textwidth]{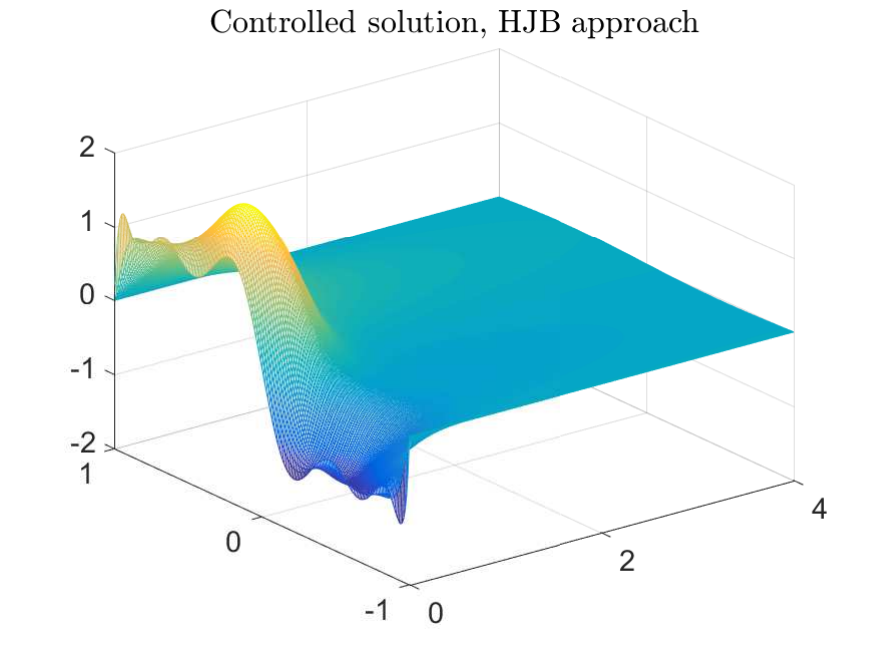}
\includegraphics[width=0.45\textwidth]{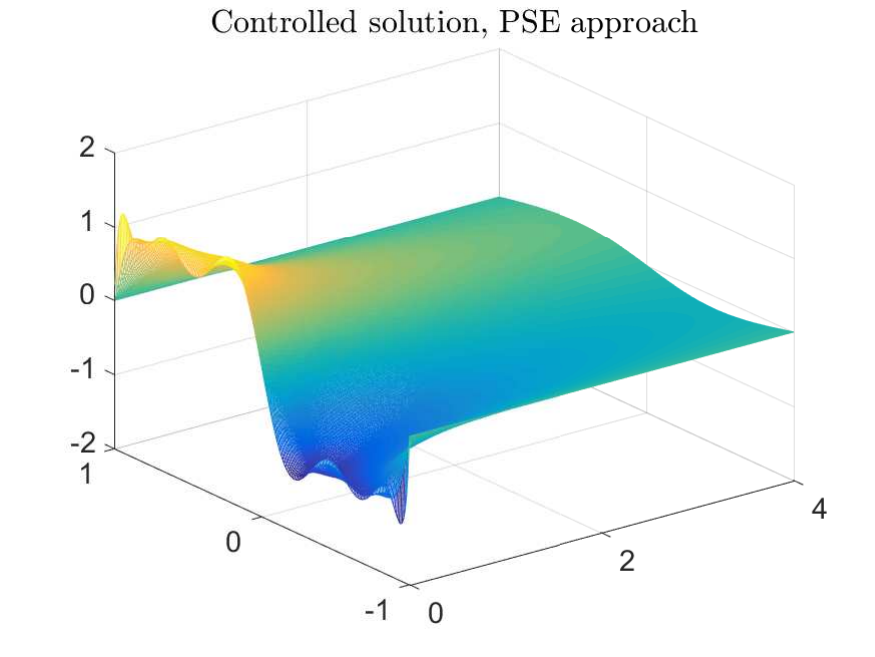}

\caption{Test 1: Viscous Burgers'-like equation. $X(\xi,0)=-\text{sign}(\xi)$. Total costs $\cJ(u,X)$: {\bf i)} Uncontrolled: $+\infty$, {\bf ii)} LQR: 7.55, {\bf iii)} PSE: 6.87, {\bf iv)} HJB: 6.25}\label{fig:test1}
\end{figure}
The feedback stabilization of Burgers' equation (without the exponential source term) has been thoroughly studied in different contexts, including the work of \cite{BK91}, and the recent work \cite{KR15}. Since our interest is the study of optimal stabilization,  we consider an additional source term $1.5X(\xi,t)e^{-0.1X(\xi,t)}$ such that the origin is not asymptotically stable. This can be appreciated in the numerical results shown in Figure \ref{fig:test1}. For this model, we consider a reduced-order state space representation of 12 states, solving a HJB equation over $\Omega=(-2,2)^{12}$. The value function is approximated with a monomial basis including both even and odd-degree polynomials up to degree 4.
In Figure \ref{fig:test1} we can compare the uncontrolled solution to the LQR- and HJB-controlled solutions, where the LQR decay is significantly slower that the one of the HJB synthesis. The HJB controller stabilizes at a higher speed, which is reflected both in the plots and in the total costs.  The HJB controller obtains a reduction of approximately 18\% with respect to the LQR cost. More importantly, the control signals differ in sign, magnitude, and speed. Such a behavior illustrates the nonlinear character of both the control problem and the feedback law.

\subsection*{Test 2: Diffusion with unstable reaction term} 
\begin{figure}[!h]
\centering
\includegraphics[width=0.45\textwidth]{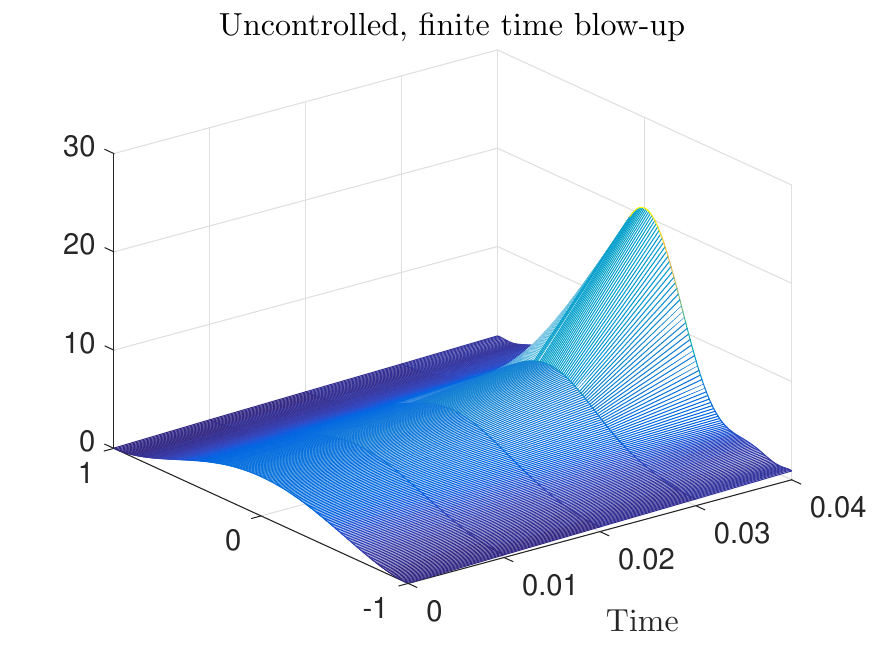}
\includegraphics[width=0.45\textwidth]{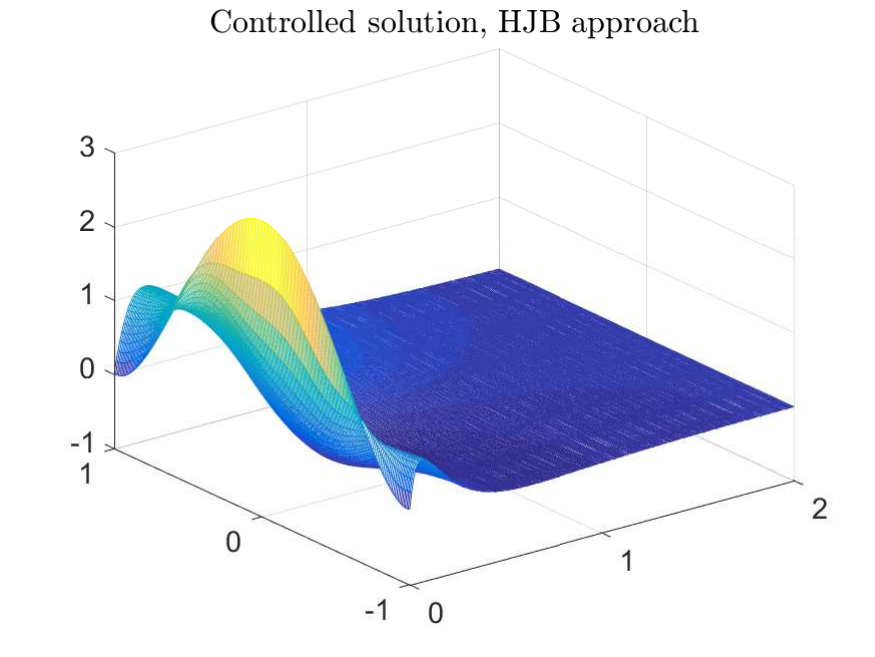}
\includegraphics[width=0.45\textwidth]{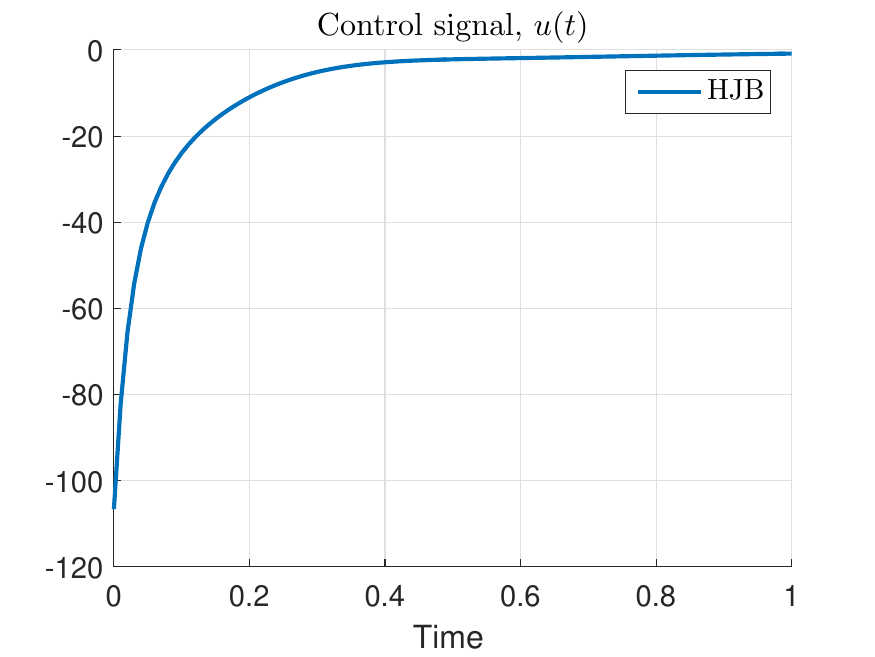}
\includegraphics[width=0.45\textwidth]{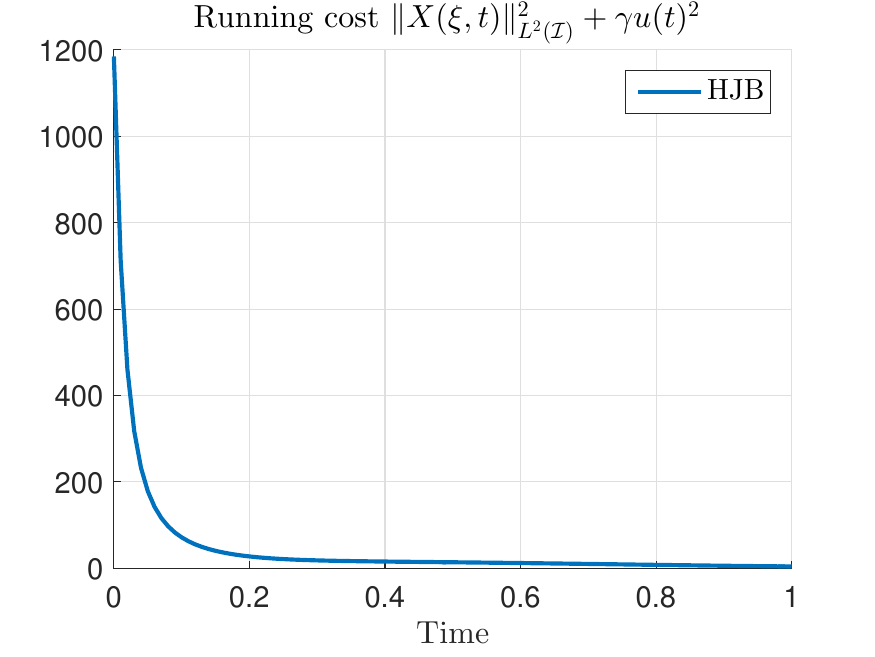}
\caption{Test 2: Diffusion with unstable reaction term. Uncontrolled dynamics leads to a finite-time blow up.}\label{fig:test2}
\end{figure}

We now turn our attention to a diffusion equation with nonlinearity $\cN(X)=X^3$ (the case with the reversed inequality sign in front of the cubic term was already treated in Subsection \ref{pdesub}),
\begin{align*}
\partial_{t}X(\xi,t) &=\sigma \partial_{\xi\xi}X(\xi,t) +X(\xi,t)^3+\chi_{\omega}(\xi)a(t)\,,\qquad\text{in}\;\cI\times\R^+\,,\\
\partial_{\xi}X(\xi_l,t)&=\partial_{\xi}X_(\xi_r,t)=0\,,\quad t\in \R^+\,,\\
X(\xi,0)&=\delta(\xi-1)^2(\xi+1)^2\,,\quad\delta\in\R^+\,, \xi \in \cI\,.
\end{align*}
We close the system with Neumann boundary conditions. The origin $X(\xi,t)\equiv 0$ is an unstable equilibrium of the uncontrolled dynamics. Any other initial condition is unstable with finite time blow-up. In this case, feedback controls can only provide local stabilization, and the purpose of this numerical test is to show that HJB-based synthesis leads to an increased closed-loop asymptotic stability region when compared to LQR, and PSE controllers. For this purpose, we compute feedback controls with the LQR, PSE and HJB approaches, for initial conditions of the form $X_0(\xi)=\delta(\xi-1)^2(\xi+1)^2$, with $\delta\in \R^+$. { The HJB feedback is computed with Algorithm \ref{alg:sga2} initialized with a nonlinear feedback control law provided by the PSE approach}. The test is carried out over $\Omega=(-2,2)^{12}$, and the value function is approximated with monomial basis elements of degree 2 and 4.
Numerical results are presented in Figure \ref{fig:test2}, for $\delta=2$  and  for a series of increased values of $\delta$ in Table \ref{tabtestcubic}. As the magnitude of the initial condition grows, the locally stabilizing LQR and PSE controllers are not able to prevent the finite blow-up of the dynamics. This eventually also happens for the HJB feedback, but at a much larger value of $\delta$ (we report the last value $\delta=4$ until which the HJB control stabilizes the dynamics).
\begin{table}[h]
\centering
$\cN(X)=X^3$, $X(\xi,0)=\delta(\xi-1)^2(\xi+1)^2$\\
\setlength{\tabcolsep}{1mm}
  \begin{tabular}{cccc}
  \hline\\
   Controller &$\delta=2$&$\delta=3$&$\delta=4$\\
   \cmidrule(lr){1-1}\cmidrule(lr){2-2}\cmidrule(lr){3-3}\cmidrule(lr){4-4}\\
   Uncontrolled & $+\infty$ & $+\infty$ & $+\infty$\\
   LQR & 4.14& $+\infty$ & $+\infty$\\
   PSE & 4.09 & 14.09 & $+\infty$ \\
   HJB & 4.06 & 13.98 & 50.36 \\
    \hline
    \\
  \end{tabular}
  \caption{Cubic source term $\cN(X)=X^3$ and increasing initial conditions. The HJB feedback law is the one which exhibits the largest closed-loop stability region. }\label{tabtestcubic}
\end{table}

\subsection*{Test 3: Newell-Whitehead equation}
\begin{figure}[!h]
\centering
\includegraphics[width=0.45\textwidth]{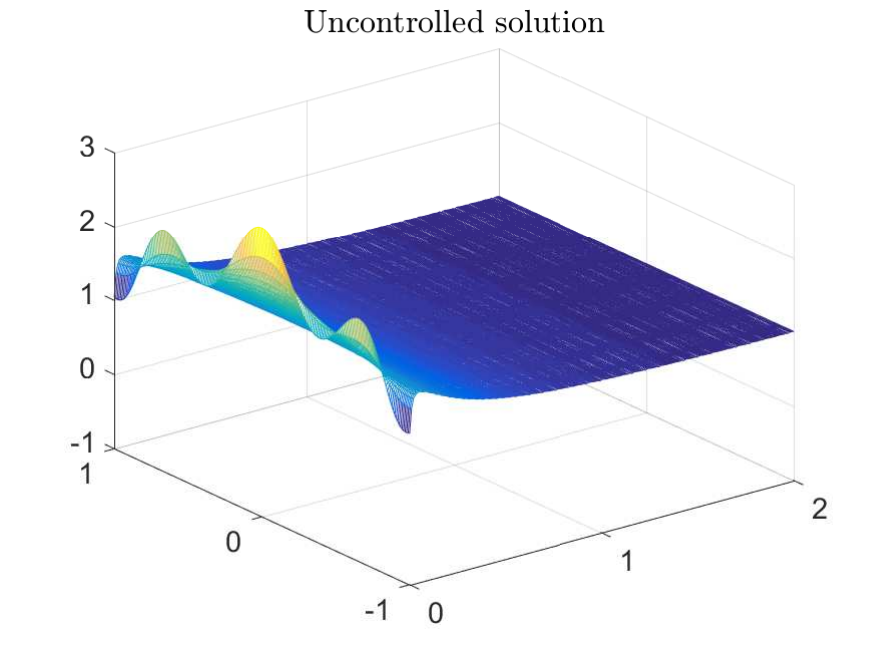}
\includegraphics[width=0.45\textwidth]{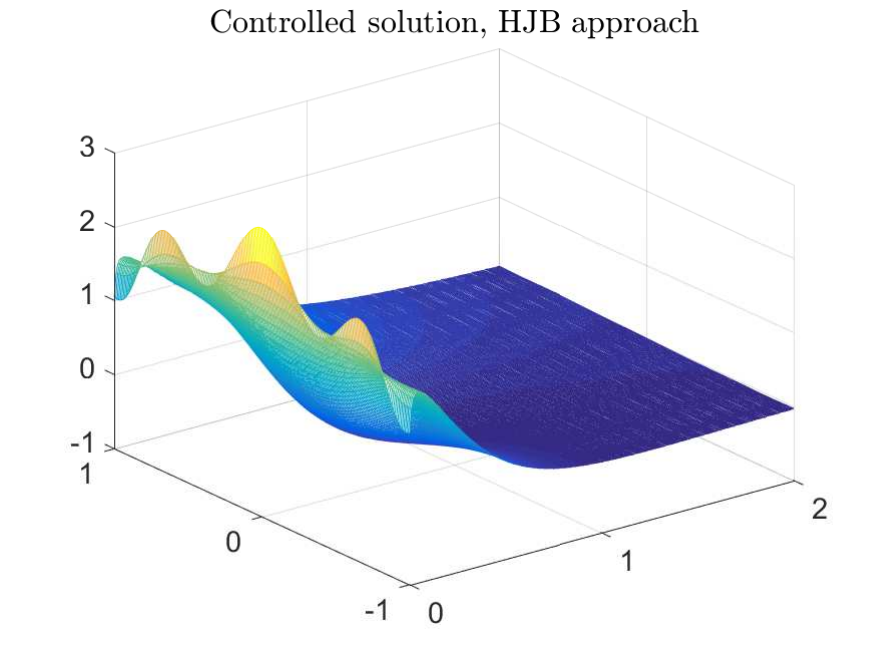}
\includegraphics[width=0.45\textwidth]{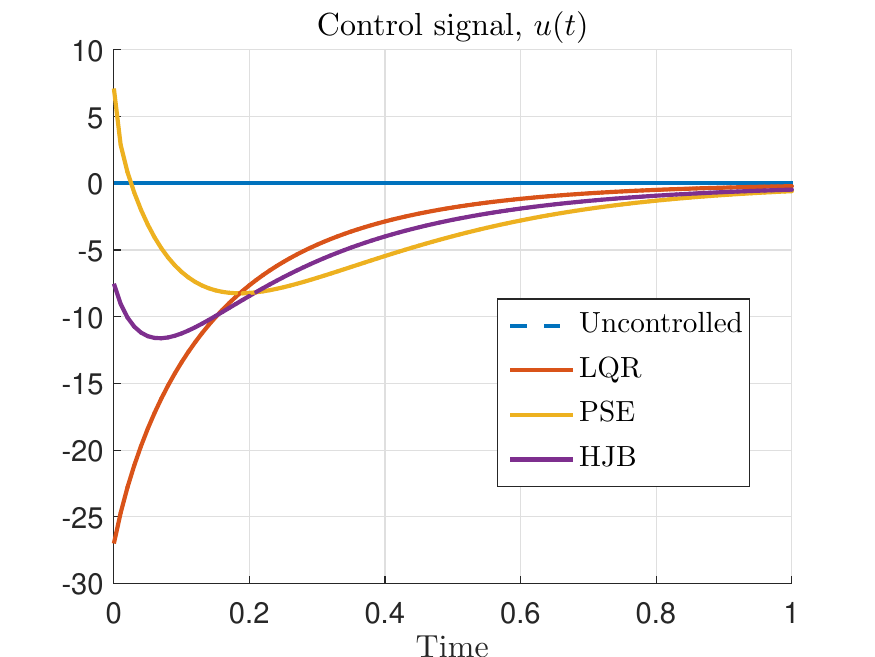}
\includegraphics[width=0.45\textwidth]{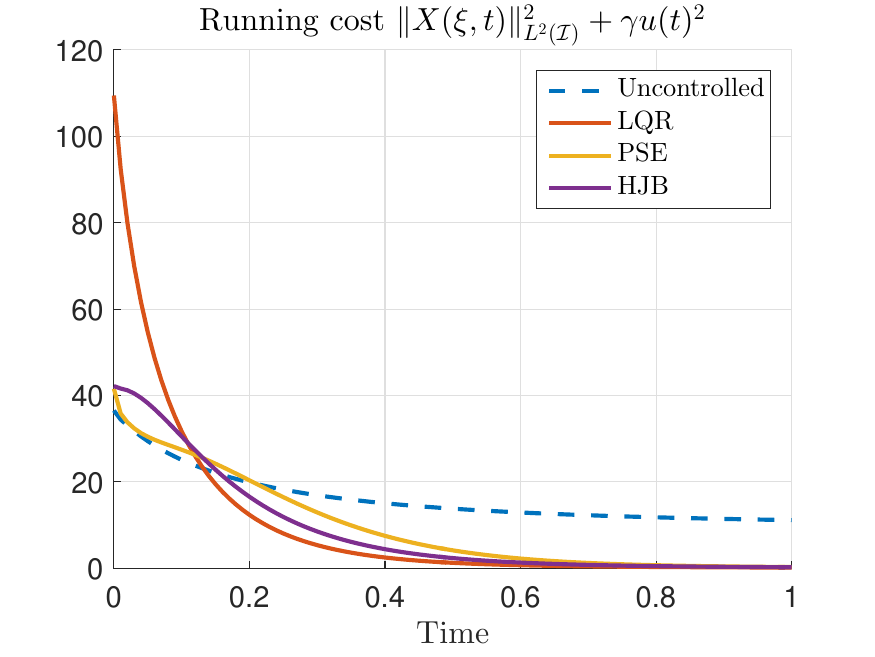}
\caption{Test 3: Newell-Whitehead equation. Initial condition $X_0(\xi)=cos(2\pi \xi)cos(\pi \xi)+2$. Uncontrolled dynamics are attracted by the stable equilibrium $X=1$. Total costs $\cJ(u,X)$  {\bf i)} Uncontrolled: $\infty$, {\bf ii)} LQR: 10.17, {\bf iii)} PSE:9.69, {\bf iv)} HJB: 8.85}\label{fig:test3}
\end{figure}

. The diffusion-reaction equation
\begin{align*}
\partial_{t}X(\xi,t) &=\sigma \partial_{\xi\xi}X(\xi,t) +X(\xi,t)(1-X(\xi,t)^2)+\chi_{\omega}(\xi)a(t)\,,\qquad\text{in}\;\cI\times\R^+\,,\\
\partial_{\xi}X(\xi_l,t)&=\partial_{\xi}X(\xi_r,t)=0\,,\quad t\in \R^+\,,\\
X(\xi,0)&=X(\xi,0)=cos(2\pi\xi)cos(\pi\xi)+\delta)\,,\quad\delta\in\R^+, \xi\in \cI\,,
\end{align*}
corresponds to a particular case of the so-called Schl\"ogl model, whose feedback stabilization has been studied in \cite{KB15,GT15}. This is a special case of a bistable system with $\pm1$ as stable and $0$ as unstable equilibria. Here we use in an essential manner that we consider Neumann boundary conditions. For Dirichlet conditions the only equilibrium is the origin. Such systems arise for instance  in Rayleigh-Benard convection and describe excitable systems such as neurons or axons. As in the previous example, the reduced state-space is chosen as $\Omega=(-2,2)^{12}$, and the basis elements for the HJB approach are even degree monomials of degree 2 and 4. Numerical results for the different controllers are shown in Figure \ref{fig:test3}. While all the feedback laws effectively stabilize the initial condition $X_0(\xi)=cos(2\pi \xi)cos(\pi \xi)+2$ to the origin, the HJB feedback has the smallest overall cost $\cJ(u,X)$. As in Test 1, it can be observed that the three feedback strategies have a considerably different transient behavior. Note that the LQR controller, which neglects the effect of the nonlinearity $\cN(X)=-X^3$, has an increased control magnitude with respect to the nonlinear controllers which are able to account the dissipative effect of the nonlinearity.
For the sake of completeness, we also consider this test case with a switch of the sign of nonlinearity, i.e., $\cN(X)=-X^3$. This case is more demanding than Test 2, as now the linear part is $\sigma \partial_{\xi\xi}X+X$. However, the performance of the controllers is similar as in Test 2, and the results are summarized in Table \ref{tabtest22}. Again, the HJB feedback law has an increased closed-loop stability region compared to the LQR and PSE controllers.
\begin{table}[h]
\centering
$\cN(X)=X^3$, $X(\xi,0)=cos(2\pi\xi)cos(\pi\xi)+\delta$\\
\setlength{\tabcolsep}{1mm}
  \begin{tabular}{cccc}
  \hline\\
   Controller &$\delta=1$&$\delta=1.5$&$\delta=2$\\
   \cmidrule(lr){1-1}\cmidrule(lr){2-2}\cmidrule(lr){3-3}\cmidrule(lr){4-4}\\
   Uncontrolled & $+\infty$ & $+\infty$ & $+\infty$\\
   LQR & 5.09& $+\infty$ & $+\infty$\\
   PSE & 4.92 & 20.02 & $+\infty$ \\
   HJB & 4.89 & 17.35 & 31.02 \\
    \hline \\
  \end{tabular}
  \caption{Test 3 with $\cN(X)=X^3$, for increasing initial conditions. Different local control strategies are not able to stabilize the dynamics for large initial conditions. The HJB control law has an increased region of the state space where it can stabilize.}\label{tabtest22}
\end{table}

\subsection*{Test 4: Degenerate Zeldovich equation} 
\begin{figure}[!h]
\centering
\includegraphics[width=0.45\textwidth]{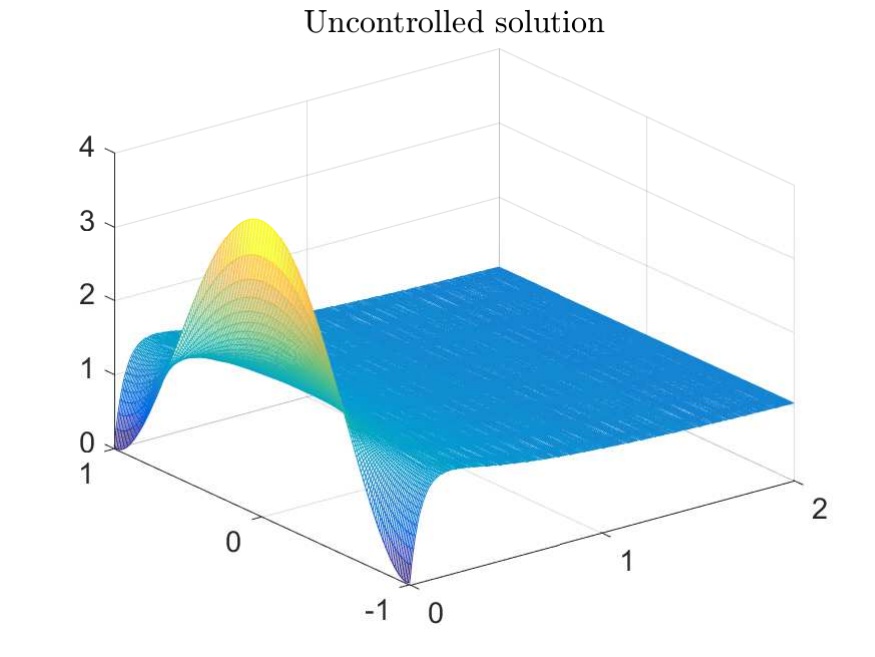}
\includegraphics[width=0.45\textwidth]{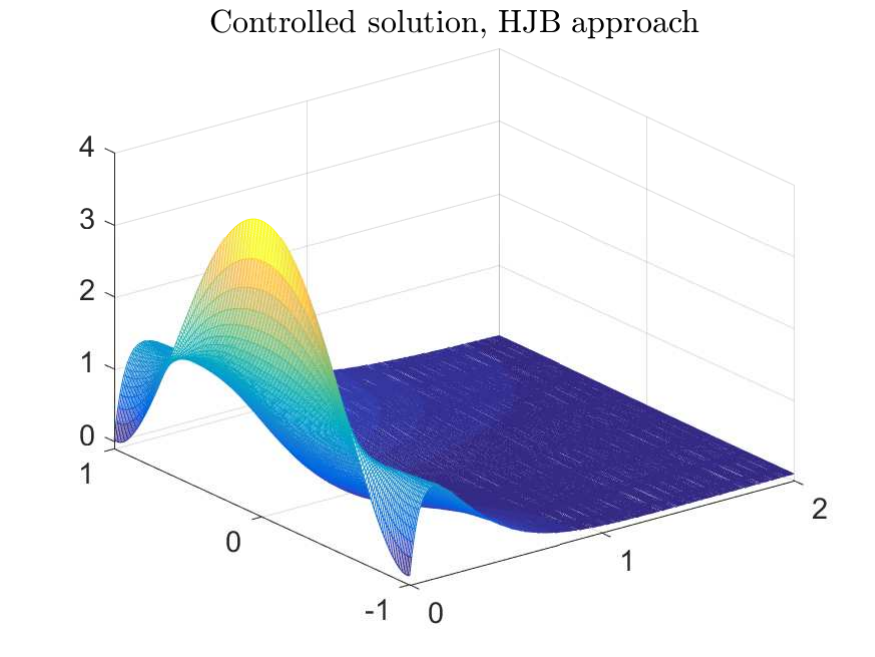}
\includegraphics[width=0.45\textwidth]{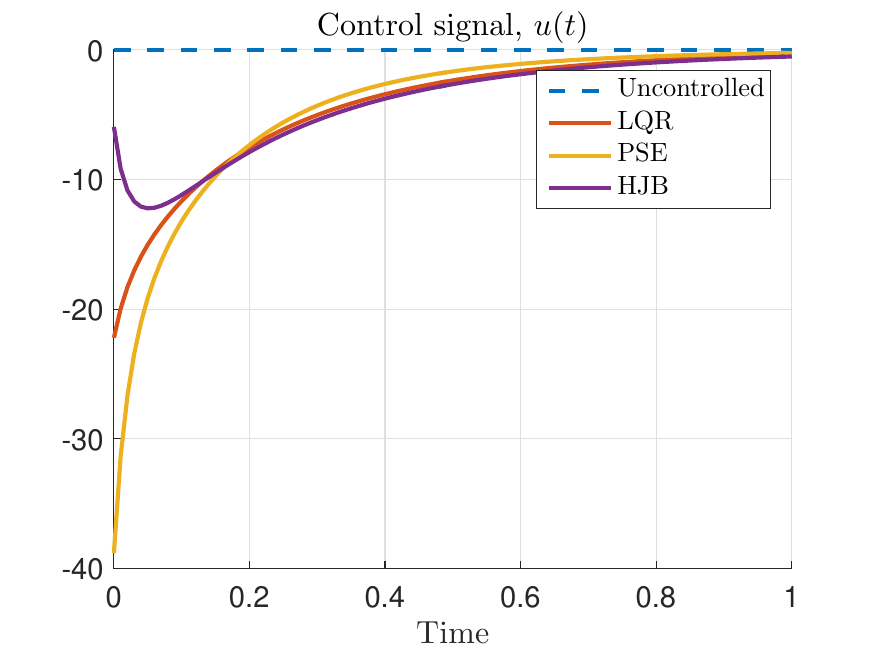}
\includegraphics[width=0.45\textwidth]{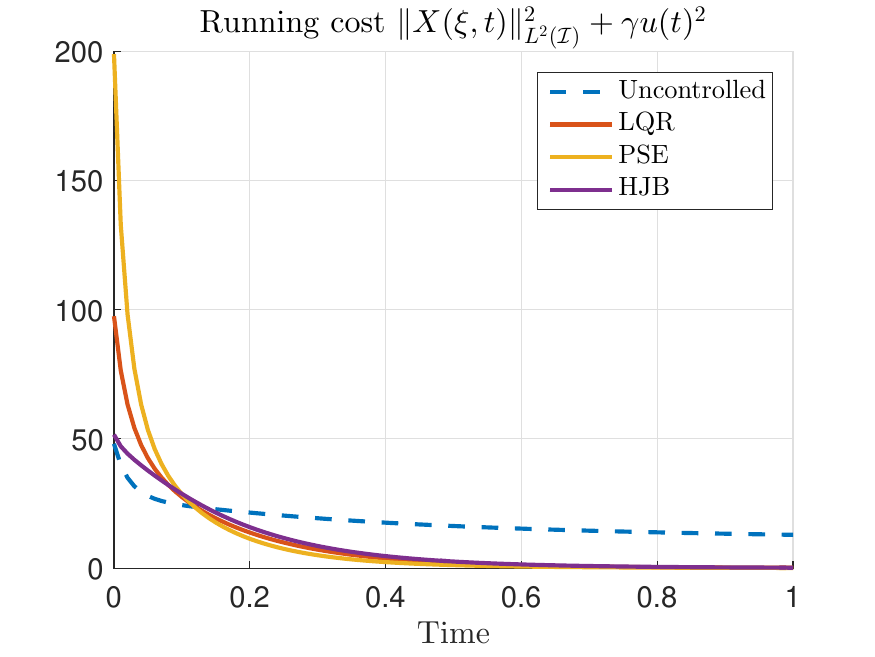}
\caption{Test 4: Degenerate Zeldovich equation. Initial condition $X_0(\xi)=4(\xi-1)^2(\xi+1)^2$. Total costs $\cJ(u,X)$ {\bf i)} Uncontrolled: $\infty$, {\bf ii)} LQR: 9.45, {\bf iii)} PSE: 11.25, {\bf iv)} HJB: 8.91 }\label{fig:test4}
\end{figure}

In this last test case, we consider the model given by
\begin{align*}
\partial_{t}X(\xi,t) &=\sigma \partial_{\xi\xi}X(\xi,t)+X(\xi,t)^2-X(\xi,t)^3+\chi_{\omega}(\xi)a(t)\,,\qquad\text{in}\;\cI\times\R^+\,,\\
\partial_{\xi}X(\xi_l,t)&=\partial_{\xi}X(\xi_r,t)=0\,,\quad t\in \R^+\,,\\
X(\xi,0)&=4(\xi-1)^2(\xi+1)^2\,,\quad \xi\in \cI\,.
\end{align*}
This equation, which arises for instance in combustion theory, has $X\equiv 1$ as stable and $X\equiv 0$ as unstable equilibria. For this case, we increase the dimension of the HJB domain to 14, i.e., $\Omega=(-2,2)^{14}$, and the basis functions are monomials of odd and even degree up to 4. Numerical results are shown in Figure \ref{fig:test4}, where it can be seen that the HJB controller yields the smaller overall cost $\cJ(u,X)$. Note that the PSE controller for this case has a diminished performance as compared even to the LQR controller. This can be explained by the fact that the PSE controller only takes into account the lowest order nonlinearity, in this case $\cN_l(X)=X^2$, neglecting the cubic term. This is a well-known drawback of this controller, and therefore justifies the need of more complex synthesis methods for nonlinear feedback design, such as the proposed HJB approach.
\section*{Concluding remarks}
A systematic technique for the computational approximation of HJB equations in optimal control problems related to semilinear parabolic equations was presented. To  partially circumvent the curse of dimensionality, the dynamics of the parabolic equation are approximated by a pseudospectral collocation  method, and the generalized HJB equation is approximated by separable multi-dimensional basis functions of a given order. The numerical results show that the feedback controls obtained by the proposed methodology differ and improve upon applying Riccati approaches to the linearized equations. The generalized HJB approach has been addressed in earlier publications,  reporting on numerical results with lower dimensions than here and in part restrained enthusiasm about the numerical performance, possibly due  to the lack of a systematic initialization procedure. For the class of problems considered in this paper the results were consistently better than Riccati approaches. The use of the discount factor path-following technique as proposed in Algorithm 2 is essential for stabilizing to unstable equilibria.

\section*{Acknowledgments}
The authors gratefully acknowledge support by the ERC advanced grant 668998
(OCLOC) under the EU's H2020 research program. Finally, D. Kalise wishes to dedicate this paper to the memory of Alexander Vasiliev.
%%%%%%%%%%%%%%%%%%%%%%%%%%%%%%%%%%%%%%%%%%%%%%%%%%%%%
% References
%%%%%%%%%%%%%%%%%%%%%%%%%%%%%%%%%%%%%%%%%%%%%%%%%%%%

\end{document}